\documentclass[reqno,11pt]{amsart}

\usepackage{amsmath,amsfonts,amssymb,amsthm,epsfig}

\usepackage{a4wide}
\usepackage{tikz-cd}
\usetikzlibrary{arrows}
\pgfarrowsdeclare{<<<}{>>>}
{
\arrowsize=0.2pt
\advance\arrowsize by .5\pgflinewidth
\pgfarrowsleftextend{-4\arrowsize-.5\pgflinewidth}
\pgfarrowsrightextend{.5\pgflinewidth}
}
{
\arrowsize=0.2pt
\advance\arrowsize by .5\pgflinewidth
\pgfpathmoveto{\pgfpointorigin}
\pgfpathlineto{\pgfpoint{-1.2mm}{-.8mm}}
\pgfusepathqstroke
\pgfpathmoveto{\pgfpointorigin}
\pgfpathlineto{\pgfpoint{-1.2mm}{.8mm}}
\pgfusepathqstroke
}

 \allowdisplaybreaks
\numberwithin{equation}{section}

\font\script=rsfs10 at 12pt
\def\eps{\varepsilon}
\def\H{{\mbox{\script H}\,\,}}
\def\comp{\subset\subset}
\def\R{\mathbb R}
\def\S{\mathbb S}
\def\A{\mathcal A}
\def\F{\mathfrak F}
\def\En{\mathfrak R}
\def\N{\mathbb N}
\def\angle#1#2#3{#1\widehat{#2}#3}
\def\sumn{\sum\nolimits}
\def\bal{\begin{aligned}}
\def\eal{\end{aligned}}
\def\proofof#1{\begin{proof}[Proof of #1]}
\def\Chi#1{\hbox{{\large $\chi$}{\Large $_{_{#1}}$}}}
\def\case#1#2{\par\noindent{\underline{\it Case~#1.}}\emph{ #2}\\}
\def\step#1#2{\par\noindent{\underline{\it Step~#1.}}\emph{ #2}\\}

\def\XXint#1#2#3{{\setbox0=\hbox{$#1{#2#3}{\int}$} \vcenter{\vspace{-1pt}\hbox{$#2#3$}}\kern-.5\wd0}}

\def\de0#1{\rule[3pt]{#1}{0.4pt} \hspace{-0.1pt} \rule[3.05pt]{0.05pt}{0.4pt} \hspace{-0.1pt} \rule[3.1pt]{0.05pt}{0.4pt} \hspace{-0.1pt} \rule[3.15pt]{0.05pt}{0.4pt} \hspace{-0.1pt} \rule[3.2pt]{0.05pt}{0.4pt} \hspace{-0.1pt} \rule[3.25pt]{0.05pt}{0.4pt} \hspace{-0.1pt} \rule[3.3pt]{0.05pt}{0.4pt} \hspace{-0.1pt} \rule[3.35pt]{0.05pt}{0.4pt} \hspace{-0.1pt} \rule[3.4pt]{0.05pt}{0.4pt} \hspace{-0.1pt} \rule[3.45pt]{0.05pt}{0.4pt} \hspace{-0.1pt} \rule[3.5pt]{0.05pt}{0.4pt} \hspace{-0.1pt} \rule[3.55pt]{0.05pt}{0.4pt} \hspace{-0.1pt} \rule[3.6pt]{0.05pt}{0.4pt} \hspace{-0.1pt} \rule[3.65pt]{0.05pt}{0.4pt} \hspace{-0.1pt} \rule[3.7pt]{0.05pt}{0.4pt} \hspace{-0.1pt} \rule[3.75pt]{0.05pt}{0.4pt} \hspace{-0.1pt} \rule[3.8pt]{0.05pt}{0.4pt} \hspace{-0.1pt} \rule[3.85pt]{0.05pt}{0.4pt} \hspace{-0.1pt} \rule[3.9pt]{0.05pt}{0.4pt} \hspace{-0.1pt} \rule[3.95pt]{0.05pt}{0.4pt} \hspace{-0.1pt} \rule[4.0pt]{0.05pt}{0.4pt} \hspace{-0.1pt} \rule[4.05pt]{0.05pt}{0.4pt} \hspace{-0.1pt} \rule[4.1pt]{0.05pt}{0.4pt} \hspace{-0.1pt} \rule[4.15pt]{0.05pt}{0.4pt} \hspace{-0.1pt} \rule[4.2pt]{0.05pt}{0.4pt} \hspace{-0.1pt} \rule[4.25pt]{0.05pt}{0.4pt} \hspace{-0.1pt} \rule[4.3pt]{0.05pt}{0.4pt} \hspace{-0.1pt} \rule[4.35pt]{0.05pt}{0.4pt} \hspace{-0.1pt} \rule[4.4pt]{0.05pt}{0.4pt} \hspace{-0.1pt} \rule[4.45pt]{0.05pt}{0.4pt} \hspace{-0.1pt} \rule[4.5pt]{0.05pt}{0.4pt} \hspace{-0.1pt} \rule[4.55pt]{0.05pt}{0.4pt} \hspace{-0.1pt} \rule[4.6pt]{0.05pt}{0.4pt} \hspace{-0.1pt} \rule[4.65pt]{0.05pt}{0.4pt} \hspace{-0.1pt} \rule[4.7pt]{0.05pt}{0.4pt} \hspace{-0.1pt} \rule[4.75pt]{0.05pt}{0.4pt} \hspace{-0.1pt} \rule[4.8pt]{0.05pt}{0.4pt} \hspace{-0.1pt} \rule[4.85pt]{0.05pt}{0.4pt} \hspace{-0.1pt} \rule[4.9pt]{0.05pt}{0.4pt} \hspace{-0.1pt} \rule[4.95pt]{0.05pt}{0.4pt} \hspace{-0.1pt} \rule[5.0pt]{0.05pt}{0.4pt} \hspace{-0.1pt} \rule[5.05pt]{0.05pt}{0.4pt} \hspace{-0.1pt} \rule[5.1pt]{0.05pt}{0.4pt} \hspace{-0.1pt} \rule[5.15pt]{0.05pt}{0.4pt} \hspace{-0.1pt} \rule[5.2pt]{0.05pt}{0.4pt} \hspace{-0.1pt} \rule[5.25pt]{0.05pt}{0.4pt} \hspace{-0.1pt} \rule[5.3pt]{0.05pt}{0.4pt} \hspace{-0.1pt} \rule[5.35pt]{0.05pt}{0.4pt} \hspace{-0.1pt} \rule[5.4pt]{0.05pt}{0.4pt} \hspace{-0.1pt} \rule[5.45pt]{0.05pt}{0.4pt} \hspace{-0.1pt} \rule[5.5pt]{0.05pt}{0.4pt} \hspace{-0.1pt} \rule[5.55pt]{0.05pt}{0.4pt} \hspace{-0.1pt} \rule[5.6pt]{0.05pt}{0.4pt} \hspace{-0.1pt} \rule[5.65pt]{0.05pt}{0.4pt} \hspace{-0.1pt} \rule[5.7pt]{0.05pt}{0.4pt} \hspace{-0.1pt} \rule[5.75pt]{0.05pt}{0.4pt} \hspace{-0.1pt} \rule[5.8pt]{0.05pt}{0.4pt} \hspace{-0.1pt} \rule[5.85pt]{0.05pt}{0.4pt} \hspace{-0.1pt} \rule[5.9pt]{0.05pt}{0.4pt} \hspace{-0.1pt} \rule[5.95pt]{0.05pt}{0.4pt} \hspace{-0.1pt} \rule[6.0pt]{0.05pt}{0.4pt}}	

\newtheorem{theorem}{Theorem}[section]
\newtheorem{lemma}[theorem]{Lemma}

\newtheorem{prop}[theorem]{Proposition}

\newtheorem{defin}[theorem]{Definition}

\begin{document}

\date{\today}

\title{Minimisers of a general Riesz-type Problem}

\author{M. Novaga}
\address{Matteo Novaga\\ Dipartimento di Matematica Universit\`a di Pisa\\ Largo Bruno Pontecorvo 5\\ 56127 Pisa\\ Italy}
\email{matteo.novaga@unipi.it}

\author{A. Pratelli}
\address{Aldo Pratelli\\ Dipartimento di Matematica Universit\`a di Pisa\\ Largo Bruno Pontecorvo 5\\ 56127 Pisa\\ Italy}
\email{aldo.pratelli@unipi.it}

\begin{abstract}
We consider sets in $\R^N$ which minimise, for fixed volume, the sum of the perimeter and a non-local term given by the double integral of a kernel $g:\R^N\setminus\{0\}\to \R^+$. We establish some general existence and regularity results for minimisers. In the two-dimensional case we show that balls are the unique minimisers in the perimeter-dominated regime, for a wide class of functions $g$.
\end{abstract}

\maketitle

\tableofcontents

\section{Introduction}

George Gamow introduced in \cite{Gam} the so-called \emph{liquid drop model} of the atomic nucleus,
corresponding to a variational problem characterized by the
competition of short-range attractive interactions, modelled by surface tension,
and long-range repulsive interactions, captured by treating the nuclear charge as
uniformly spread throughout the nucleus.
A competition of the attractive forces which try
to minimise the interfacial area of the nucleus, and the repulsive Coulombic forces that try to
spread the charges apart makes the nucleus unstable at sufficiently large atomic numbers,
resulting in nuclear fission. On the other hand, in Gamow's model the nucleus is stable and spherical for small atomic numbers.
We refer to \cite{CMT}, and references therein, for a comprehensive introduction to this model.

In this paper we are interested in a generalisation of Gamow's model and we consider the energy 
\[
\F(\Omega) := P(\Omega) + \iint_{\Omega\times\Omega} g(y-x)\,dy\,dx
\]
for a set $\Omega\subseteq\R^N$, where 
$g:\R^N\setminus \{0\}\to \R^+$ is a function modelling the repulsing interaction of the set $\Omega$ with itself.
The Gamow's model corresponds to the choice $N=3$ and $g(x)=1/|x|$.
When $g$ is a radially decreasing function, there is a clear competition between the two terms of the energy: On one side, the perimeter term favours concentration, and it is minimised by the ball; on the other side the minimisation of the repulsive term favours disgregation of the set. 

We observe that the problem is not invariant by rescaling. 
Indeed, for $m>0$ the functional $\F$ satisfies the scaling property
\begin{eqnarray*}
\F(m \Omega) &=& P(m\Omega) + \iint_{m\Omega\times m\Omega} g(y-x)\,dy\,dx
\\
&=& m^{N-1} \left(P(\Omega) + m^{N+1} \iint_{\Omega\times\Omega} g\left(m(y-x)\right)\,dy\,dx\right).
\end{eqnarray*}
In particular, if we assume that $g(mx)=m^\alpha [g(x) + o(1)]$ as $m\to 0$, with $\alpha>-(N+1)$,
minimising $\F$ among sets of volume $m\ll 1$ is closely related (and completely equivalent if $g(x)=|x|^\alpha$) 
to minimising the functional 
\[
\F_\eps(\Omega) := P(\Omega) + \eps \iint_{\Omega\times\Omega} g(y-x)\,dy\,dx
\]
for $\eps = m^{N+1+\alpha} \ll 1$, among sets $\Omega$ of fixed volume. 
In this paper we will focus on the minimisation of the functional $\F_\eps$.

A remarkable fact is that, in some cases, volume-constrained minimisers of $\F_\eps$ are actually balls for $\eps$ small enough. 
This has been recently established in \cite{KM1,KM2} (see also~\cite{CS,FFMMM,MS,Jul,FL})
for the physically relevant case of a negative power, that is, $g(x)=|x|^{\alpha}$ for some $\alpha\in (-N,0)$.
Here we show that the same is true for a wide class of functions $g$, in the $2$-dimensional case.
More precisely, we shall prove the following result:

\begin{theorem}\label{2dcase}
Let $g:\R^2\setminus\{0\}\to\R^+$ be a radial, decreasing, and positive definite function (in the sense of Definition~\ref{defpos}) such that
\begin{equation}\label{ipLip}
\int_0^1 g(t)\,dt<+\infty\,.
\end{equation}
Then there exists some $\bar\eps>0$ such that, for every $0<\eps<\bar\eps$, the unique minimiser of $\F_\eps$ in the class
\[
\A := \Big\{ \Omega\subseteq\R^N\text{ measurable}:\, |\Omega|=\omega_N\Big\}
\]
is the unit ball.
\end{theorem}
\noindent We notice that assumption \eqref{ipLip} implies that the functional is finite on sets of finite perimeter, and in particular it is satisfied by $g(x)=|x|^\alpha$ for $\alpha \in (-1,0)$.

\smallskip 

As shown in \cite{KM1,KM2} for the liquid drop model, existence fails for $\eps$ big enough, since minimisers tend to split
in two or more components (nuclear fission) which then move far apart one from the other in oder to decrease the nonlocal energy.
To capture this phenomenon and describe the shape of the components, 
it is convenient to introduce a generalised energy defined as

\begin{equation}\label{genprob}
\widetilde\F_\eps(\Omega) := \inf_{H\in\N} \widetilde \F_{\eps,H}(\Omega)\,,
\end{equation}
where
\[
\widetilde\F_{\eps,H}(\Omega) :=\inf \left\{ \sum_{i=1}^H \F_\eps(\Omega^i):\, \Omega=\bigcup_{i=1}^H \,\Omega^i,\, 
\Omega^i\cap\Omega^j=\emptyset\quad \text{for}\, 1\leq i\neq j\leq H\right\}\,.
\]
Notice that in this functional the interaction between different components is not evaluated, which corresponds to consider them ``at infinite distance''. By considering $\widetilde\F_\eps$ instead of $\F_\eps$, we can prove the following general existence result:

\begin{prop}[Existence of generalised minimisers]\label{genex}
Let $g$ be an admissible decreasing function (in the sense of Definitions~\ref{def:admissible} and~\ref{defdec}). 
For every $\eps>0$ there exists a minimiser of $\widetilde\F_\eps$ in the class $\A$. 
More precisely, there exist a set $E\in\A$ and a partition $E=\cup_{i=1}^H E^i$, with pairwise disjoint sets $E^i$, such that
\[
\widetilde \F_\eps(E) = \sumn_{h=1}^H \F_\eps(E^h) = \inf \Big\{ \widetilde \F_\eps(\Omega):\, \Omega\in\A\Big\}\,.
\]
Moreover, for each $1\leq \bar i \leq H$ the set $E^{\bar i}$ is a minimiser of both the standard and the generalised energy for its volume, i.e. it satisfies
\begin{equation}\label{eachEiminim}
\widetilde\F_\eps(E^{\bar i}) = \F_\eps(E^{\bar i}) = \min \Big\{\widetilde \F_\eps(\Omega):\, \Omega\subseteq\R^N:\, |\Omega|=|E^{\bar i}|\Big\}\,.
\end{equation}
\end{prop}

We observe that the existence result in Proposition \ref{genex} was previously obtained only in some particular cases, see for instance
 \cite{KMNa} (see also \cite{KMN}) for $N=3$ and $g(x)=1/|x|$, and \cite{MS} for $N=2$ and $g(x)=\chi_{[\delta,\infty)}(|x|)/|x|^3$, with $\delta>0$.
 
\smallskip

The plan of the paper is the following: in Section \ref{secnot} we give some preliminary definitions and prove some sufficient conditions for a function $g$ to be positive definite. In Section \ref{sec:prelim} we show existence and regularity of generalised minimisers, thus proving Proposition \ref{genex}. Finally, in Section \ref{secballs} we prove Theorem \ref{2dcase} on the minimality of balls for small values of the parameter $\eps$, in two dimensions.

\smallskip 

{\bf Acknowledgements.} The authors acknowledge partial support from the INDAM-GNAMPA and from the PRIN Projects 2017TEXA3H and 2017BTM7SN.
 
\section{Notation and preliminary results}\label{secnot}

In the following we denote by $B_r$ the ball of radius $r$ and centre in the origin, and we set for simplicity $B:=B_1$. Given two measurable sets $F,\,G\subseteq\R^N$, we also let
\begin{align*}
\En(F,G):=\int_F \int_G g(y-x)\,dy\,dx\,, && \En(F) := \En(F,F)\,.
\end{align*}

\begin{defin}[Admissible functions]\label{def:admissible}
A function $g:\R^N\setminus\{0\}\to\R^+$ is \emph{admissible} if $\En(B)<+\infty$. 
\end{defin}

\noindent Notice that a radial function $g$ is admissible if and only if (see also Lemma~\ref{lemadm})
\[
\int_0^1 g(t) t^{N-1}\,dt<+\infty\,.
\]

\begin{defin}[Decreasing function]\label{defdec}
We say that the admissible function $g:\R^N\setminus\{0\}\to\R^+$ is \emph{decreasing} if for every $x\in\R^N\setminus\{0\}$ and every $\lambda>1$ one has
\[
g(\lambda x)\leq \lambda g(x)\,.
\]
\end{defin}

\begin{lemma}
For every admissible $g$ there exists a constant $C$ such that
\[
\En(B_r) \leq C r^N
\]
for every $0<r\leq 1$.
\end{lemma}

\begin{proof}
Let us denote by $Q_\ell$ the cube centred at the origin, with half-side $\ell$. 
Since for any $m\in\N$ the cube $Q_1$ is the essentially disjoint union of $m^N$ cubes of side $1/m$, of course
\[
\En(Q_1) \geq m^N \En(Q_{1/m})\,.
\]
For any $0<r\leq 1$, let us call $m$ the integer part of $1/r$, so that $(2r)^{-1}<m\leq r^{-1}$. The above estimate then gives
\[
\En(Q_1) \geq m^N \En(Q_{1/m}) \geq \frac 1{(2r)^N} \En(Q_r)\,,
\]
thus
\[
\En(Q_r) \leq 2^N \En(Q_1) r^N\,,
\]
from which the thesis follows since $B_r \subseteq Q_r$.
\end{proof}

We now recall the concept of positive definite functions (see for instance~\cite[Chapter 11]{Rudin}).

\begin{defin}[Positive definite functions]\label{defpos}
A function $g\in L^1_{\rm loc}(\R^N)$ is \emph{positive definite} if
\[
\iint_{\R^N\times\R^N} g(y-x)f(x)f(y)\,dy\,dx \ge 0 \qquad \forall f\in C_c(\R^N).
\]
Equivalently, $g$ is positive definite if
\begin{equation}\label{ineqposit}
\En(F) + \En(G) \geq 2\En(F,G)\,,
\end{equation}
for every two bounded measurable sets $F,\,G\subseteq\R^N$.
\end{defin}

\begin{lemma}\label{subhpos}
Every admissible, subharmonic function $g:\R^N\setminus\{0\}\to \R^+$ such that
\begin{equation}\label{infinf}
\lim_{|x|\to\infty} g(x) = \inf_{x\in\R^N} g(x) < +\infty
\end{equation}
is positive definite.
\end{lemma}

\begin{proof}
By $L^1$ continuity of the integral, a simple approximation argument shows that it is sufficient to show, for every positive constants $m_1,\,m_2,\, L,\,R$, that the inequality~(\ref{ineqposit}) is valid for every couple of sets of finite perimeter $F,\,G\subseteq \R^N$ with
\begin{align}\label{HypBan}
{\rm diam}(F)\leq R\,, && {\rm diam}(G)\leq R\,, &&
|F| \leq m_1\,, && |G|\leq m_2\,, &&
P(F)+P(G)\leq L\,.
\end{align}
Let us introduce
\[
K:= \inf \Big\{\En(F)+\En(G)-2\En(F,G),\, F,\,G \subseteq\R^N \hbox{ satisfy~(\ref{HypBan})} \Big\}\,.
\]
If $K\geq 0$, there is nothing to prove, we can then assume by contradiction that $K<0$. In this case, the assumption~(\ref{infinf}) immediately implies that $K$ is actually a minimum.\par

Let then $F,\,G$ be two sets realizing the minimum. We reduce ourselves to the case that $F\cap G=\emptyset$, since otherwise we can call $\widetilde F=F\setminus G$ and $\widetilde G=G\setminus F$, and observe that also the sets $\widetilde F,\,\widetilde G$ satisfy~(\ref{HypBan}), and
\[
\En(\widetilde F)+\En(\widetilde G)-2\En(\widetilde F,\widetilde G)=\En(F)+\En(G)-2\En(F,G) = K\,.
\]
For every $v\in\R^N$, we call $F_v=v+F$, and we define $f:\R^N\to\R$ as
\[
f(v) := \En(F_v)+\En(G)-2\En(F_v,G)= \En(F)+\En(G)-2\En(F_v,G)\,.
\]
By approximation, we can assume that $g$ is strictly subharmonic. As a consequence, the function $f$ is strictly superharmonic at $v=0$ (here we use that $F\cap G=\emptyset$), against the fact that $0$ is a minimum of $f$ by construction. The contradiction shows the thesis.
\end{proof}

\noindent Notice that $g(x)=|x|^{-\alpha}$, with $0<\alpha<N-1$, satisfies the assumptions of Lemma ~\ref{subhpos}.

\begin{lemma}\label{lemBoc}
If $g:\R^N\setminus \{0\}\to \R^+$ is admissible, satisfies
\[
\limsup_{|x|\to\infty} g(x) < +\infty,
\] 
and
its Fourier transform $\hat g$ is a nonnegative Borel measure, then $g$ is positive definite.
\end{lemma}

\begin{proof}
Let $F,G\subset\R^N$ be bounded measurable sets, and let $f := \Chi{F}-\Chi{G}\in L^\infty(\R^N)$.
Notice tha $\hat f$ is continuous on $\R^N$. Then, by Plancherel's Theorem we get that
\[
\En(F)+\En(G)-2\En(F,G) = \langle g*f,f\rangle_{L^2(\R^N)} = \int_{\R^N} \hat f(\xi)^2 \, d\hat g(\xi)\geq 0,
\]
which gives \eqref{ineqposit}.
\end{proof}

\noindent Notice that $g(x)=e^{-\kappa |x|^2}|x|^{-\alpha}$, with $\kappa\ge 0$ and 
$0<\alpha<N$, satisfies the assumptions of Lemma ~\ref{lemBoc}.

\section{Existence and regularity of generalised minimisers}\label{sec:prelim}

In this section we collect some general properties of the functional $\F_\eps$. A first, easy fact is the following.

\begin{lemma}\label{omog}
Let $\Omega\subseteq\R^N$ be a set of finite perimeter, and let $\lambda>1$. If $g$ is decreasing (in the sense of Definition~\ref{defdec}), then
\[
\F_\eps(\lambda \Omega) \leq \lambda^{2N} \F_\eps(\Omega)\,.
\]
\end{lemma}
\begin{proof}
This is immediate because
\[
P(\lambda\Omega) = \lambda^{N-1} P(\Omega) \leq \lambda^{2N} P(\Omega)
\]
and, since $g$ is decreasing, we have
\[\begin{split}
\En(\lambda\Omega) &= \iint_{(\lambda \Omega)^2} g(y-x)\,dy\,dx = \lambda^{2N} \iint_{\Omega\times\Omega} g(\lambda(y-x))\,dy\,dx
\leq \lambda^{2N} \iint_{\Omega\times\Omega} g(y-x)\,dy\,dx \\
&= \lambda^{2N} \En(\Omega)\,.
\end{split}\]
\end{proof}

A consequence of the above estimate is the next geometric lemma, which allows to ``cut and paste'' an excessively long and thin set decreasing its energy.

\begin{lemma}\label{cap}
For every $\overline m\in\R$ there exists a positive constant $L>0$ such that the following holds. Let $E\subseteq\R^N$, and let $a<b$ be two numbers with $b>a+2L$ and such that
\[
\Big| \Big\{x\in E:\, a\leq x_1\leq b\Big\}\Big| \leq \overline m\,.
\]
There exist then two numbers $a< a^+< a+L$ and $b-L< b^-< b$ such that, by calling $E^-= E\setminus [a^+,b^-]\times\R^{N-1}$ and $m= |E|\setminus |E^-|\leq \overline m$, one has
\begin{equation}\label{cutandpaste}
\F_\eps(E^-) \leq \F_\eps(E) - \frac 12\, N\omega_N^{1/N} m^{\frac{N-1}N}\,.
\end{equation}
\end{lemma}
\begin{proof}
Let us fix $\overline m$, and let $L=L(N,\overline m)$ be a constant, to be precised later. Let also $E,\,a,\,b$ be as in the claim. For almost every $t\in\R$, we set
\[
\sigma(t) = \H^{N-1}\Big(E \cap \big\{ x\in\R^N:\, x_1=t\big\}\Big)\,.
\]
Let then $c=(a+b)/2$, and let us call
\[
\varphi(t) = \int_t^c \sigma(s)\,ds\,,
\]
so that $\varphi(c)=0$ and $\varphi(a)\leq \overline m$. We claim the existence of some $a< a^+< a+L<c$ such that
\begin{equation}\label{choicea+}
\sigma(a^+) \leq \frac 18\, N\omega_N^{1/N} \Big(\Big| \Big\{x\in E:\, a^+<x_1<c\Big\}\Big| \Big)^{\frac{N-1}N}=: \frac 18\, N\omega_N^{1/N} m_1^{\frac{N-1}N}\,.
\end{equation}
Indeed, for every $t\in (a,a+L)$, if the above inequality is false with the choice $a^+=t$ then
\[
-\varphi'(t)=\sigma(t)>\frac 18\, N\omega_N^{1/N} \varphi(t)^{\frac{N-1}N}\,.
\]
Since the power $\tfrac{N-1}N$ is strictly less than $1$, if a positive, decreasing function $\varphi:[a,d]\to \R^+$ satisfies
\[
\left\{\begin{array}{c}
\varphi(a)\leq \overline m\,,\\
|\varphi'(t)| > \frac 18\, N\omega_N^{1/N} \varphi(t)^{\frac{N-1}N}\,,
\end{array}\right.
\]
then the length $d$ is bounded by a constant, depending only on $\overline m$ and $N$. We can then choose any number $L$ larger than this constant, and the existence of some $a< a^+< a+L$ satisfying~(\ref{choicea+}) is guaranteed.\par

Similarly, we have the existence of some $b-L< b^-<b$ such that
\[
\sigma(b^-) \leq \frac 18\, N\omega_N^{1/N} \Big(\Big| \Big\{x\in E:\, c<x_1<b^-\Big\}\Big| \Big)^{\frac{N-1}N}=:\frac 18\, N\omega_N^{1/N} m_2^{\frac{N-1}N}\,.
\]
We call then $E^-=E\setminus [a^+,b^-]\times\R^{N-1}$ and, to conclude, we have to establish~(\ref{cutandpaste}). Notice that, calling $F=E\setminus E^-$, we have by construction
\[
|F|=m=m_1+m_2\,,
\]
so that by the isoperimetric inequality
\[
P(F) \geq N\omega_N^{1/N} m^{\frac{N-1}N} \,.
\]
We conclude then
\[\begin{split}
P(E^-)&\leq P(E) -P(F) + 2\big(\sigma(a^+)+\sigma(b^-)\big)
\leq P(E) -P(F) + \frac 14\, N\omega_N^{1/N} \Big(m_1^{\frac{N-1}N}+m_2^{\frac{N-1}N}\Big)\\
&\leq P(E) -P(F) + \frac 12\, N\omega_N^{1/N} \big(m_1+m_2\big)^{\frac{N-1}N}
\leq P(E) -\frac 12\, N\omega_N^{1/N} \big(m_1+m_2\big)^{\frac{N-1}N}\,.
\end{split}\]
Since of course $\En(E^-) \leq \En(E)$ because $E^-\subseteq E$, we deduce~(\ref{cutandpaste}) and the proof is concluded.
\end{proof}

As a consequence, we obtain the following uniform boundedness result.

\begin{lemma}\label{unifbdd}
For every $m\in\R$ there exist two constants $R>0$ and $\overline H\in\N$, depending only on $N,\,m$ and $\eps$, such that if $g$ is decreasing then
\[
\inf \Big\{ \F_\eps(\Omega):\, \Omega\subseteq\R^N,\,|\Omega|=m\Big\} \geq \inf \Big\{ \widetilde\F_{\eps,\overline H}^R(\Omega):\, \Omega\subseteq\R^N,\,|\Omega|=m\Big\}\,,
\]
where
\[
\widetilde\F_{\eps,\overline H}^R(\Omega)=
\inf \bigg\{ \sumn_{i=1}^{\overline H} \F_\eps(\Omega^i):\, \Omega=\cup_{i=1}^{\overline H} \,\Omega^i,\, \Omega^i\cap\Omega^j=\emptyset,\,
 {\rm diam}\, \Omega^i\leq R\, \ \forall\, 1\leq i\neq j\leq \overline H \bigg\}\,.
\]
\end{lemma}

\begin{proof}
Let $M\in\N$ be a natural number, only depending on $N,\,m$ and $\eps$ and to be specified later, and let us call $\overline m=m/M$. Let $E \subseteq\R^N$ be any bounded set with $|E|=m$ such that
\begin{equation}\label{omalop}
\F_\eps(E) \leq \inf \Big\{ \F_\eps(\Omega):\, \Omega\in\A\Big\} + \frac{N\omega_N^{1/N}}3\, \bigg(\frac m{M^2}\bigg)^{\frac{N-1}N}\,,
\end{equation}
which is possible because clearly the infimum is reached by a sequence of bounded sets. Let $t_0 < t_1<t_2<\, \cdots\,< t_{M-1}< t_M$ be real numbers such that
\[
\Big| E \cap (t_i,t_{i+1})\times\R^{N-1}\Big| =\overline m\,,
\]
for every $0\leq i\leq M-1$. Let $L=L(\overline m,N)$ be given by Lemma~\ref{cap}.\par

For every $0\leq i\leq M-1$, if $t_{i+1}-t_i\leq 2L$ then we set $I_i=\emptyset$, while otherwise we call $a=t_i,\, b=t_{i+1}$, we let $a^+$ and $b^-$ be given by Lemma~\ref{cap}, and we set $I_i=[a^+,b^-]$. Let then $m_i = \big| E \cap I_i\times\R^{N-1}\big|\leq \overline m$. We claim that
\begin{equation}\label{mismall}
m_i \leq \frac m{M^2}\,.
\end{equation}
This is clearly true if $I_i=\emptyset$. Otherwise, we let
\begin{align*}
E' = \lambda \Big(E\setminus \big(I_i\times \R^{N-1}\big)\Big)\,, && \hbox{where} && \lambda= \bigg(\frac m{m-m_i}\bigg)^{1/N}\,.
\end{align*}
Notice now that $m_i/m\leq 1/M$ by construction. Hence, provided that $M$ is large enough, by Lemma~\ref{omog} and~(\ref{cutandpaste}) we have
\[\begin{split}
\F_\eps(E') &\leq \bigg(\frac m{m-m_i}\bigg)^2 \F_\eps\Big(E\setminus \big(I_i\times \R^{N-1}\big)\Big)
\leq \bigg(1+3\,\frac{m_i}m\bigg) \bigg(\F_\eps(E) - \frac 12\, N\omega_N^{1/N} m_i^{\frac{N-1}N}\bigg)\\
&\leq \F_\eps(E) - \frac 13\, N\omega_N^{1/N} m_i^{\frac{N-1}N}\,.
\end{split}\]
Since by construction we have that $|E'|=m$, this estimate, together with~(\ref{omalop}), immediately implies~(\ref{mismall}).\par

Let us now call
\[
\widetilde E = E \setminus \Big(\bigcup\nolimits_{i=0}^{M-1} I_i \times \R^{N-1}\Big)\,,
\]
and let $\mu=\sum_{i=0}^{M-1} m_i$, so that $|\widetilde E|=m-\mu$. We can apply Lemma~\ref{cap} $M$ times, finding
\[
\F_\eps(\widetilde E) \leq \F_\eps(E) - \frac 12\, N\omega_N^{1/N} \sumn_{i=0}^{M-1} m_i^{\frac{N-1}N}
\leq \F_\eps(E) - \frac 12\, N\omega_N^{1/N} \mu^{\frac{N-1}N}\,.
\]
The set $F=(m/(m-\mu))^{1/N} \widetilde E$ has then volume $m$, and by Lemma~\ref{omog} we obtain
\[
\F_\eps(F)\leq \bigg(\frac m{m-\mu}\bigg)^2\, \F_\eps(\widetilde E)
\leq \bigg(\frac m{m-\mu}\bigg)^2\, \bigg(\F_\eps(E) - \frac 12\, N\omega_N^{1/N} \mu^{\frac{N-1}N}\bigg)\,.
\]
We deduce $\F_\eps(F)\leq \F_\eps (E)$ as soon as $\mu$ is small enough, which in turn is true if $M$ is big enough, since by~(\ref{mismall}) we get $\mu\leq m/M$.\par

Notice that, by construction, $\widetilde E$ is the union of at most $M+1$ parts, and by Lemma~\ref{cap} each of them has horizontal width at most equal to $2ML$. As a consequence, recalling again that $\mu\leq m/M$, $F$ is the union of at most $M+1$ parts, each of them with horizontal width at most equal to $3ML$.\par

Repeating in the obvious way the same argument $N-1$ times, in order to get boundedness of the pieces in all the $N$ directions, we clearly obtain the existence of two constants $R>0$ and $\overline H\in\N$, and of a set $G\subseteq\R^N$ with $|G|=m$ and $\F_\eps(G)\leq \F_\eps(E)$, so that $G$ is the disjoint union of sets $G_i,\, 1\leq i \leq \overline H$, and each $G_i$ has diameter at most $R$. Since of course
\[
\F_\eps(G) \geq \sum_{i=1}^{\overline H} \F_\eps(G_i) \geq \widetilde\F_{\eps,\overline H}^R(G)\,,
\]
the proof is concluded.
\end{proof}

We can now prove Proposition \ref{genex}.

\proofof{Proposition~\ref{genex}}
We fix $\eps>0$ and we split the proof in few steps.

\step{I}{Reduction to the case $H(n)=H'$.}
First of all, we claim the existence of a natural number $H'$ and of a sequence $\{G_n\}_{n\in\N}\subseteq\A$ such that
\begin{equation}\label{goalstep1}
\inf \Big\{ \widetilde \F_\eps(\Omega):\, \Omega\in\A\Big\} = \lim_{n\to \infty} \widetilde\F_{\eps,H'}(G_n)\,.
\end{equation}
To do so, we let $H'$ be a large integer, only depending on $N$ and $\eps$ and to be specified later, and we start with a generic minimising sequence, that is, a sequence of sets $\{\Omega_n\}_{n\in\N}\subseteq\A$ such that
\begin{equation}\label{optimalsequence}
K:=\inf \Big\{ \widetilde \F_\eps(\Omega):\, \Omega\in\A\Big\} = \lim_{n\to\infty} \widetilde\F_\eps(\Omega_n)\,.
\end{equation}
For every $n\in\N$, we can take a number $H(n)\in\N$ and a subdivision $\Omega_n=\Omega_n^1\cup \Omega_n^2\cup \cdots \cup \Omega_n^{H(n)}$ so that
\begin{equation}\label{almostopt}
\widetilde\F_\eps(\Omega_n) > \bigg(1-\frac 1{n+1}\bigg) \sumn_{i=1}^{H(n)} \F_\eps(\Omega_n^i)\,.
\end{equation}

Let us now focus on a given $n\in\N$. For brevity of notation, we call $H=H(n)$ and $m_i = \big| \Omega_n^i\big|$ for every $1\leq i\leq H$, and we assume without loss of generality that $m_i$ is decreasing with respect to $i$. We observe that by~(\ref{almostopt})
\[\begin{split}
\widetilde\F_\eps(\Omega_n) &\geq \frac 12\, \sumn_{i=1}^H P(\Omega_n^i) \geq \frac 12\, \sumn_{i=1}^H N\omega_N^{1/N} m_i^{\frac{N-1}N}
\geq \frac 1{2\sqrt[N]{m_1}}\, \sumn_{i=1}^H N\omega_N^{1/N} m_i\\
&=\frac 1{2\sqrt[N]{m_1}}\,N\omega_N^{\frac{N+1}N}\,,
\end{split}\]
which by~(\ref{optimalsequence}) implies
\begin{equation}\label{boundonm1}
m_1 \geq \bigg(\frac{N\omega_N^{\frac{N+1}N}}{4K}\bigg)^N
\end{equation}
for every $n$ large enough. For every such $n$, we define then the set
\begin{align*}
G_n = \bigcup\nolimits_{i=1}^{H'} \lambda \Omega_n^i\,, && \hbox{being} && \lambda = \bigg(\frac{\omega_N}{\omega_N-\sum_{i>H'} m_i}\bigg)^{1/N}\leq 1 + C_1 \sumn_{i=H'+1}^H m_i\,,
\end{align*}
where $C_1$ is a constant, depending only on $N$ and $K$, so actually on $N$ and $\eps$, whose existence is guaranteed by~(\ref{boundonm1}). Notice that, by the definition of $\lambda$, also the set $G_n$ belongs to $\A$. By Lemma~\ref{omog}, we can estimate
\begin{equation}\label{fromabove}
\widetilde\F_{\eps,H'}(G_n)\leq \sum_{i=1}^{H'} \F_\eps(\lambda\Omega_n^i) \leq \lambda^{2N} \sum_{i=1}^{H'} \F_\eps(\Omega_n^i)
\leq \Big(1+C_2 \sum_{i=H'+1}^H m_i\Big)\sum_{i=1}^{H'} \F_\eps(\Omega_n^i)\,,
\end{equation}
where $C_2=C_2(N,\eps)$ is another constant. On the other hand, by~(\ref{almostopt}) we have
\[\begin{split}
\widetilde\F_\eps(\Omega_n) &> \bigg(1-\frac 1{n+1}\bigg) \sumn_{i=1}^H \F_\eps(\Omega_n^i)
\geq \bigg(1-\frac 1{n+1}\bigg) \sumn_{i=1}^{H'} \F_\eps(\Omega_n^i) + \frac 12\sumn_{i=H'+1}^H P(\Omega_n^i)\\
&\geq \bigg(1-\frac 1{n+1}\bigg) \sumn_{i=1}^{H'} \F_\eps(\Omega_n^i) + N\omega_N^{1/N} \sumn_{i=H'+1}^H m_i^{\frac{N-1}N}\,.
\end{split}\]
Putting this estimate together with~(\ref{fromabove}), and keeping in mind that, by~(\ref{optimalsequence}) and~(\ref{almostopt}), for $n$ large enough we surely have
\[
\sumn_{i=1}^{H'} \F_\eps(\Omega_n^i) \leq 2K\,,
\]
we get
\begin{equation}\label{rightest}
\widetilde\F_{\eps,H'}(G_n)-\widetilde\F_{\eps,H}(\Omega_n) \leq 2K \Big(C_2 \sumn_{i=H'+1}^H m_i + \frac 1{n+1}\Big) - N\omega_N^{1/N} \sumn_{i=H'+1}^H m_i^{\frac{N-1}N}\,.
\end{equation}
We are now in position to define $H'$ as an integer so large that
\[
H' \geq \bigg(\frac{2KC_2}{N}\bigg)^N\,.
\]
Notice that $H'$ only depends on $N$ and on $\eps$. With this choice of $H'$, and observing that $m_i\leq \omega_N/H'$ for every $i>H'$ by construction, from~(\ref{rightest}) we deduce
\[
\widetilde\F_{\eps,H'}(G_n)\leq \widetilde\F_{\eps,H}(\Omega_n) + \frac{2K}{n+1}
\]
for every $n$ large enough. Keeping in mind~(\ref{optimalsequence}), the sequence $\{G_n\}_{n\in\N}\subseteq\A$ satisfies~(\ref{goalstep1}), hence this step is concluded.

\step{II}{Bound on $\F_\eps$ by means of $\widetilde\F_\eps$.}
In this step we show that for every $m>0$ there exist an integer $\overline H\in\N$, a bounded set $E$ with $|E|=m$, and a subdivision $E=\bigcup_{j=1}^{\overline H} E^j$ such that
\begin{equation}\label{goalstep2}
\widetilde\F_\eps(E)\leq \sumn_{j=1}^{\overline H} \F_\eps(E^j)\leq\inf \Big\{ \F_\eps(\Omega):\, \Omega\subseteq\R^N,\, |\Omega|=m\Big\}\,.
\end{equation}
To do so, we let $R$ and $\overline H$ be given by Lemma~\ref{unifbdd}, and we let $\{\Omega_n\}_{n\in\N}$ be a sequence of sets of volume $m$ such that
\begin{equation}\label{rc1}
\inf \Big\{ \F_\eps(\Omega):\, \Omega\subseteq\R^N,\, |\Omega|=m\Big\} \geq \lim_{n\to\infty} \widetilde\F_{\eps,\overline H}^R(\Omega_n)\,,
\end{equation}
where $\widetilde\F_{\eps,\overline H}^R$ is defined as in Lemma~\ref{unifbdd}. 

For every $n\in\N$, we can then write $\Omega_n=\Omega_n^1\cup \Omega_n^2\cup \cdots \cup \Omega_n^{\overline H}$ where the sets $\Omega_n^i$ are all disjoint for $1\leq i\leq \overline H$, ${\rm diam}(\Omega_n^i)\leq R$, and
\begin{equation}\label{rc2}
\sumn_{i=1}^{\overline H} \F_\eps(\Omega_n^i) \leq \widetilde \F_{\eps,\overline H}^R (\Omega_n) + \frac 1 n\,.
\end{equation}
By compactness, up to a subsequence we have constants $m_i,\, 1\leq i \leq \overline H$, so that
\begin{align*}
m_i = \lim_{n\to \infty} |\Omega_n^i|\quad \forall\, 1\leq i \leq \overline H\,, && m=\sumn_{i=1}^{\overline H} m_i\,.
\end{align*}

Let us now fix $1\leq j\leq \overline H$, and let us concentrate on the sets $\{\Omega_n^j\}_{n\in\N}$. Since they all have diameter less than $R$, up to translations we can assume that they are all contained in a fixed ball with radius $R$. The characteristic functions $f_n=\Chi{\Omega_n^j}$ have then uniformly bounded supports, and they are bounded in $BV$, since the perimeter of every set $\Omega_n^j$ is clearly less than $\F_\eps(\Omega_n^j)$. Up to a subsequence, we can then assume that $f_n$ weakly* converge in $BV$ to some function $f$. Since the convergence is in particular strong in $L^1$, then also $f$ is the characteristic function of a bounded set with volume $m_j$, that we call $E^j$. By the lower-semicontinuity of the perimeter under weak* 
$BV$-convergence, and the continuity of $\En$ under strong $L^1$ convergence, we obtain that
\begin{equation}\label{rc3}
\F_\eps(E^j) \leq \liminf_{n\to\infty} \F_\eps(\Omega_n^j)\,.
\end{equation}
Since the diameters of the sets $E^j$ are all bounded, up to a translation we can assume that the sets are disjoint, and we can call $E=\cup_{j=1}^{\overline H} E^j$, which is a finite union of bounded sets, so bounded itself. By construction, $|E|=m$, and by~(\ref{rc3}), (\ref{rc2}) and~(\ref{rc1}) we have
\[\begin{split}
\widetilde\F_\eps(E) &\leq \sumn_{j=1}^{\overline H} \F_\eps(E^j) 
\leq \sumn_{j=1}^{\overline H} \liminf_{n\to\infty} \F_\eps(\Omega_n^j)
\leq \liminf_{n\to\infty} \sumn_{j=1}^{\overline H} \F_\eps(\Omega_n^j)
\leq \liminf_{n\to\infty} \widetilde \F_{\eps,\overline H}^R (\Omega_n)\\
&\leq \inf \Big\{ \F_\eps(\Omega):\, \Omega\subseteq\R^N,\, |\Omega|=m\Big\}\,,
\end{split}\]
so~(\ref{goalstep2}) is proved and this step is concluded.

\step{III}{Proof of the existence.}
Thanks to Step~I, we have a sequence of sets $\{G_n\}_{n\in\N}$, and a subdivision of each of the sets as $G_n=G_n^1\cup G_n^2\cup \cdots \cup G_n^{H'}$, so that
\begin{equation}\label{optimality}
\inf \Big\{ \widetilde \F_\eps(\Omega):\, \Omega\in\A\Big\} = \lim_{n\to \infty} \sumn_{i=1}^{H'} \F_\eps(G_n^i)\,.
\end{equation}
As in Step~II, up to a subsequence there exist constants $\mu_i,\, 1\leq i \leq H'$, such that
\begin{align*}
\mu_i = \lim_{n\to \infty} |G_n^i|\quad \forall\, 1\leq i \leq H'\,, && \omega_N=\sumn_{i=1}^{H'} \mu_i\,.
\end{align*}
If we now define
\[
K_i := \inf \Big\{ \F_\eps(\Omega):\, \Omega\subseteq\R^N,\, |\Omega|=\mu_i\Big\}\,,
\]
we have
\begin{equation}\label{sumKi}
\inf \Big\{ \widetilde \F_\eps(\Omega):\, \Omega\in\A\Big\} = \sumn_{i=1}^{H'} K_i\,,
\end{equation}
being one inequality trivial, and the other one a consequence of~(\ref{optimality}).\par

By Step~II, for every $1\leq i\leq H'$ there exist an integer $\overline H(i)\in\N$, a bounded set $E_i\subseteq\R^N$ with $|E_i|=\mu_i$, and a subdivision $E_i=\bigcup_{j=1}^{\overline H(i)} E_{i,j}$ in pairwise disjoint sets so that~(\ref{goalstep2}) holds with $m=\mu_i$, that is,
\begin{equation}\label{eachopt}
\widetilde\F_\eps(E_i)\leq \sumn_{j=1}^{\overline H(i)} \F_\eps(E_{i,j}) \leq K_i\,.
\end{equation}
Since the sets $E_i$ are bounded, up to translations we can assume them to be disjoint, so that the set $E=\cup_{i=1}^{H'} E_i$ has volume $\omega_N$. By construction, the set $E$ is the disjoint union of all the sets $E_{i,j}$ with $1\leq i\leq H'$ and $1\leq j\leq\overline H(i)$, let us call all these sets $E^h$ with $1\leq h\leq H$ and $H=\sum_{i=1}^{H'} \overline H(i)$. By~(\ref{sumKi}) and~(\ref{eachopt}) we have then
\[
\widetilde\F_\eps(E) \leq \sumn_{h=1}^H \F_\eps(E^h)\leq \sumn_{i=1}^{H'} K_i = \inf \Big\{ \widetilde \F_\eps(\Omega):\, \Omega\in\A\Big\}\,.
\]
hence the set $E$ is a minimiser of $\widetilde\F_\eps$ and the subdivision $E=\bigcup_{h=1}^H E^h$ is optimal, so the proof of the existence is concluded.

\step{IV}{Proof of~(\ref{eachEiminim})}
To conclude the proof, we only have to establish the validity of~(\ref{eachEiminim}) for a given $1\leq \bar i \leq H$. The first equality is obvious, since a non-trivial subdivision of $E^{\bar i}$ having less generalised energy than the whole $E^{\bar i}$ could be used to build another subdivision of $E$ strictly better than the optimal one, which is impossible. Concerning the second equality, if it was false then there would be a bounded set $F^{\bar i}$ with the same volume as $E^{\bar i}$ and strictly less generalised energy, say
\begin{equation}\label{hic}
\widetilde\F_\eps( F^{\bar i}) = \widetilde\F_\eps( E^{\bar i}) - \eta
\end{equation}
with some $\eta>0$. Moreover, for each $i\neq \bar i$ it is possible to find a bounded set $F^i$ with the same volume as $E^i$ and
\begin{equation}\label{hoc}
\widetilde\F_\eps( F^i) < \widetilde\F_\eps( E^i) + \frac \eta H\,.
\end{equation}
Since all the sets $F^i$ and $F^{\bar i}$ are bounded, up to a translation they are disjoint, so that the set $F=\cup_{j=1}^H F^j$ belongs to $\A$. Moreover, by~(\ref{hic}) and~(\ref{hoc}) we would have
\[
\widetilde\F_\eps (F) \leq \sum_{j=1}^H \widetilde\F_\eps (F^j) < \sum_{j=1}^H \widetilde\F_\eps(E^j) =\widetilde\F_\eps(E)\,,
\]
against the optimality of $E$.
\end{proof}

A simple observation is that the infima of $\F_\eps$ and of $\widetilde\F_\eps$ are equal if $g$ is vanishing at infinite.

\begin{lemma}\label{g0eq}
If $\lim_{|x|\to \infty} g(x)=0$, then
\begin{equation}\label{inf=}
\inf \Big\{ \widetilde\F_\eps(\Omega):\, \Omega\in\A\Big\}= \inf \Big\{ \F_\eps(\Omega):\, \Omega\in\A\Big\}\,.
\end{equation}
\end{lemma}
\begin{proof}
Since for every set $\Omega$ one has $\widetilde\F_\eps(\Omega)\leq \F_\eps(\Omega)$, one inequality is emptily true without any assumption on $g$. Concerning the other inequality, let $\Omega\in\A$ be any bounded set. By definition, for every $\delta>0$ there exists a subdivision $\Omega=\Omega_1\cup\Omega_2 \cup\cdots\cup\Omega_H$ such that
\begin{equation}\label{goodsub}
\sumn_{i=1}^H \F_\eps(\Omega_i) \leq \widetilde\F_\eps(\Omega)+\delta\,.
\end{equation}
Since $g$ is vanishing at infinite, and since the sets $\Omega_i$ are bounded, for every $1\leq i\leq H$ we can define $\Omega_i'$ as a suitable translation of $\Omega_i$ in such a way that the sets $\Omega_i'$ are pairwise disjoint and
\[
g(y-x)<\delta\qquad \forall\, x\in \Omega'_i,\, y\in\Omega'_j,\, i\neq j\,.
\]
Setting then $\Omega'=\cup_{i=1}^H \Omega_i'$, we have
\[
\En(\Omega') \leq \sumn_{i=1}^H \En(\Omega_i') + \omega_N^2 \delta
=\sumn_{i=1}^H \En(\Omega_i) + \omega_N^2 \delta\,.
\]
Also by~(\ref{goodsub}) we derive
\[\begin{split}
\F_\eps(\Omega') &= P(\Omega') +\eps\En(\Omega') \leq \eps\omega_N^2\delta + \sumn_{i=1}^H P(\Omega_i) + \eps\En(\Omega_i) 
= \eps\omega_N^2\delta +\sumn_{i=1}^H \F_\eps(\Omega_i) \\
&\leq \widetilde\F_\eps(\Omega) + (\eps\omega_N^2+1)\delta\,.
\end{split}\]
Since of course the infimum of $\widetilde\F_\eps$ can be reached by a sequence of bounded sets, the thesis is concluded.
\end{proof}

Keep in mind that the vanishing assumption on $g$ is true for all the most interesting functions $g$, as well as for the physically relevant ones. Observe that Lemma~\ref{g0eq} is done without the assumption on $g$ to be decreasing, hence we cannot apply Proposition~\ref{genex}. This is why the two terms in~(\ref{inf=}) are both infima, not necessarily minima. On the other hand, if $g$ is both decreasing and vanishing, then the result of Lemma~\ref{unifbdd} is in fact an equality, not just an inequality, and Step~III in the proof of Proposition~\ref{genex} is not needed.

We present now a simple but useful fact.

\begin{lemma}\label{FvG}
If $g$ is decreasing, then there exists a continuous and increasing function $\varphi:\R^+\to\R^+$, with $\varphi(0)=0$, such that for every two sets $F,\,G\subseteq\R^N$ one has
\begin{equation}\label{thFvG}
\En(F,G) \leq |F| \varphi(|G|)\,.
\end{equation}
\end{lemma}
\begin{proof}
For every positive number $t\in\R^+$, let us call
\[
E(t) = \Big\{x\in\R^N:\, g(x)\geq t\Big\}\cup \{0\}\,.
\]
Since $g$ is decreasing, this set is star-shaped in $0$. By construction, we have the inclusion $E(t)\subseteq E(s)$ whenever $t\geq s$, and moreover $|E(t)|\to 0$ if $t\to +\infty$ because $g$ is real-valued on $\R^N\setminus \{0\}$ and decreasing. As a consequence, for every $\sigma>0$ we can find a Borel set $D(\sigma)$ such that
\begin{align}\label{lpg}
|D(\sigma)| = \sigma\,, && g(x)\geq g(y)\quad \forall\, x\in D(\sigma),\, y\notin D(\sigma)\,.
\end{align}
We can also assume that the sets $D(\sigma)$ are ordered by inclusion, that is, $D(\sigma)\subseteq D(\sigma')$ whenever $\sigma\leq \sigma'$. We can then define $\varphi:(0,+\infty)\to [0,+\infty]$ as
\[
\varphi(\sigma) = \int_{D(\sigma)} g(y)\,dy\,.
\]
By construction, the function $\varphi$ is continuous and increasing. Moreover, by~(\ref{lpg}), for every set $G\subseteq\R^N$ with $|G|=\sigma$ we have
\[
\int_G g(y)\, dy \leq \int_{D(\sigma)} g(y)\,dy= \varphi(\sigma)\,.
\]
As a consequence, for every $x\in\R^N$ and every set $G\subseteq\R^N$, we have
\[
\int_{y\in G} g(y-x)\,dy = \int_{y\in x+G} g(y)\,dy \leq \varphi(|G|)\,,
\]
thus
\[
\En(F,G) = \int_{x\in F} \int_{y\in G} g(y-x)\,dy\, dx \leq |F| \varphi(|G|)\,,
\]
that is, we obtained~(\ref{thFvG}). To conclude, we only have to check that $\varphi$ is real-valued, and that $\varphi(t)\to 0$ if $t\to 0$.\par

To do so, let us notice that $\min\big\{ g(x),\, x\in \overline{B_{1/2}}\setminus\{0\}\big\}$ is defined and strictly positive since $g$ is l.s.c., strictly positive and decreasing. Keeping in mind that $|E(t)|\to 0$ if $t\to +\infty$ and recalling~(\ref{lpg}), we deduce that $D(\bar\sigma)\subseteq B_{1/2}$ for some small $\bar\sigma$. As a consequence,
\[\begin{split}
\En(B) &= \int_{x\in B} \int_{y\in B} g(y-x)\,dy\,dx \geq \int_{x\in B_{1/2}} \int_{y\in B_{1/2}} g(y)\,dy\,dx 
\geq \frac{\omega_N}{2^N}\, \int_{y \in D(\bar\sigma)} g(y)\,dy\\
&= \frac{\omega_N}{2^N}\, \varphi(\bar\sigma)\,,
\end{split}\]
which implies that $\varphi(\bar\sigma)<+\infty$ since $g$ is admissible in the sense of Definition~\ref{def:admissible}. By construction, we deduce then that $\varphi$ is real-valued. Finally, since the sets $D(\sigma)$ are ordered by inclusion, we have that
\[
\varphi(\sigma) = \int_{D(\sigma)} g(y)\,dy \to 0\qquad \hbox{for $\sigma\to 0$}\,,
\]
hence the proof is concluded.
\end{proof}

We can conclude the section by showing the regularity of (generalised) minimisers.

\begin{prop}[Regularity of minimisers]
Assume that $g$ is decreasing in the sense of Definition~\ref{defdec}, and that $E\in\A$ is a minimiser of the generalised problem~(\ref{genprob}) for some $\eps>0$. Then $E$ is a set of class ${\rm C}^{1,\frac 12-\delta}$ for every $\delta>0$.
\end{prop}

\begin{proof}
By Proposition~\ref{genex}, we know that $E$ is a finite, disjoint union of volume-constrained minimisers of the standard energy $\F_\eps$. 
As a consequence, it is enough to show that minimisers of $\F_\eps$ are regular.\par

Let then $\eps>0$ be given, and let $E\subseteq\R^N$ be a set such that
\begin{equation}\label{probE}
\F_\eps(E) = \min \Big\{ \F_\eps(\Omega):\, \Omega\subseteq\R^N,\, |\Omega|=|E|\Big\}\,.
\end{equation}
We have to show that $E$ is ${\rm C}^{1,\frac 12-\delta}$ regular for every $\delta>0$. By standard regularity theory (see for instance~\cite{Tam,DavSem}), it is enough to show that $E$ is $\omega$-minimal with $\omega(r)=r^{1-\delta}$ for every $\delta>0$. In other words, we have to show the existence of some $\bar r>0$ such that, for every ball $B(x,r)\subseteq\R^N$ with radius $r<\bar r$ and for every set $F\subseteq\R^N$ with $F\Delta E \comp B(x,r)$, one has
\[
P(E)\leq P(F) + r^{N-\delta}\,.
\]
Let us then assume by contradiction the existence of a ball $B(x,r)$ and of a set $F\subseteq\R^N$ such that $F\Delta E\comp B(x,r)$ and
\begin{equation}\label{regu1}
P(E) > P(F) + r^{N-\delta}\,.
\end{equation}
Let us call $\alpha=|E|-|F|$, and notice that of course $|\alpha| \leq \omega_N r^N$. 
By classical properties of sets of finite perimeter (see for instance~\cite{Tam}), it is possible to find a set $G\subseteq\R^N$ such that
\begin{align}\label{regu2}
G\Delta E\comp \R^N\setminus \overline{B(x,r)}\,, && |G|=|E|+\alpha\,, && |G\Delta E|\leq 2|\alpha|\,, && P(G) \leq P(E) + C |\alpha|\,,
\end{align}
for a geometric constant $C=C(N)$. We can then define the competitor $\Omega=\big(F\cap B(x,r)\big) \cup \big(G\setminus B(x,r)\big)$. By construction, we have that $|\Omega|=|E|$, hence $\Omega$ is an admissible set for the minimisation problem~(\ref{probE}). Moreover, by~(\ref{regu1}) and~(\ref{regu2}) we have
\begin{equation}\label{estper}\begin{split}
P(\Omega)-P(E) &= P(F) - P(E) + P(G)-P(E)
\leq C|\alpha|-r^{N-\delta}
\leq C\omega_N r^N-r^{N-\delta}\\
&\leq -\frac 12\, r^{N-\delta}\,,
\end{split}\end{equation}
where the last inequality is true as soon as $\bar r$ has been chosen small enough, only depending on $N$ and on $\delta$.\par

Let us define $E^+=\Omega\setminus E$ and $E^-=E\setminus \Omega$, so that $\Omega = (E \cup E^+) \setminus E^-$ and $E^+\cap E=\emptyset$, while $E^-\subseteq E$. By construction and~(\ref{regu2}), we have
\[
|E^+| + |E^-| \leq |B(x,r)| + 2|\alpha| \leq 3\omega_N r^N\,.
\]
As a consequence, keeping in mind Lemma~\ref{FvG}, we have
\[
\En(\Omega)-\En(E) = \En(E,E^+) - \En(E,E^-)+\En(\Omega,E^+)-\En(\Omega,E^-)
\leq 6\omega_Nr^N \varphi(|E|)\,.
\]
Putting together this estimate with~(\ref{estper}) we get
\[
\F_\eps(\Omega) = P(\Omega) +\eps \En(\Omega) \leq \F_\eps(E) +6\eps\omega_Nr^N \varphi(|E|)-\frac 12\, r^{N-\delta} < \F_\eps(E)\,,
\]
where the last inequality again holds true as soon as $\bar r$ has been chosen small enough, only depending on $N$, $\delta$, $\eps$ and $\varphi(|E|)$. Since this estimate is against the optimality of $E$ in~(\ref{probE}), we have found the searched contradiction and the proof is concluded.
\end{proof}

\section{Properties of minimisers in the perimeter-dominated regime}\label{secballs}

We now consider the situation in which the coefficient $\eps$ in the energy $\F_\eps$ is small. We can see that, still without the radial assumption on $g$, if $\eps\to 0$ then the minimising sets converge to the unit ball in the $L^1$ sense. Moreover, they have good topological properties, in particular they are connected and without holes. 

In the $2$- dimensional case we will prove Theorem~\ref{2dcase}, in which we show that, adding the radial assumption and a couple of technical assumptions, the solution is \emph{exactly} the ball for $\eps\ll 1$.\par

We start by showing the connectedness of solutions.

\begin{lemma}[Solutions are connected for $\eps\ll 1$]\label{solconn}
If $g$ is decreasing in the sense of Definition~\ref{defdec}, then there exists $\eps_1>0$ such that, for every $0<\eps<\eps_1$, every minimiser $\Omega$ of $\widetilde\F_\eps$ is connected, or equivalently, if $E\subseteq\Omega$ is such that $0<|E|<\Omega$, then
\begin{equation}\label{connectedness}
P(\Omega) < P(E) + P(\Omega\setminus E)\,.
\end{equation}
In addition, one has $\F_\eps(\Omega)=\widetilde\F_\eps(\Omega)$, hence $\Omega$ is also a minimiser of $\F_\eps$ and $\min\F_\eps=\min\widetilde\F_\eps$.
\end{lemma}

Before presenting the proof of Lemma \ref{solconn}, let us observe that, for a decreasing function $g$, Proposition~\ref{genex} ensures the existence of minimisers for $\widetilde\F_\eps$, while existence of minimisers of $\F_\eps$ is in general false. The present lemma shows that, for $\eps\ll 1$, any minimiser of $\widetilde\F_\eps$ is actually also a minimiser of $\F_\eps$ (so in particular minimisers of $\F_\eps$ exist), and moreover $\min \{\widetilde\F_\eps\}=\min\{\F_\eps\}$. In particular this shows that, if $\eps$ is small enough, equality~(\ref{inf=}) holds even without the assumption that $g$ is vanishing, so even if Lemma~\ref{g0eq} cannot be applied.
\proofof{Lemma~\ref{solconn}}
Let $\eps>0$ be given, and let $\Omega\in\A$ be any minimiser of $\widetilde\F_\eps$. By Proposition~\ref{genex}, which ensures the existence of a minimiser, we also obtain the existence of an optimal subdivision, that is, we can write $\Omega=\bigcup_{i=1}^H \Omega^i$ so that the sets $\Omega^i$ are pairwise disjoint, and
\begin{equation}\label{optsubdiv}
\widetilde\F_\eps(\Omega) = \sumn_{i=1}^H \F_\eps(\Omega^i)\,.
\end{equation}
We have then
\[
N\omega_N + \eps \En(B) = \F_\eps(B) \geq \widetilde \F_\eps(B) \geq \widetilde\F_\eps(\Omega)\geq \sumn_{i=1}^H P(\Omega^i) \geq N\omega_N^{1/N} \sumn_{i=1}^H |\Omega^i|^{\frac{N-1}N}\,,
\]
which implies
\begin{equation}\label{wrongconv}
\sumn_{i=1}^H |\Omega^i|^{\frac{N-1}N} \leq \bigg(\sumn_{i=1}^H |\Omega^i|\bigg)^{\frac{N-1}N} + \frac\eps{N\omega_N^{1/N}}\, \En(B)\,.
\end{equation}
Let $\overline m$ be a very small constant, depending on $N$ and on $\eps$, to be specified later. Recalling that $\En(B)$ is a number and that $s\mapsto s^{\frac{N-1}N}$ is strictly concave, there exists $\eps_1$ such that, as soon as $\eps<\eps_1$, the estimate~(\ref{wrongconv}) implies that (up to renumbering the sets $\Omega^i$)
\[
m:=\omega_N- \big|\Omega^1\big| <\overline m\,.
\]
Let us define
\[
E = \bigg(\frac{\omega_N}{\omega_N-m}\bigg)^{1/N}\, \Omega^1\,,
\]
which belongs to $\A$ by construction. Applying Lemma~\ref{omog} and recalling that $m<\overline m$, as soon as $\overline m$ is small enough we have then
\[
\F_\eps(E) \leq \bigg(\frac{\omega_N}{\omega_N-m}\bigg)^2 \F_\eps(\Omega^1) \leq \F_\eps(\Omega^1) + \frac{3m}{\omega_N}\,\F_\eps(\Omega^1)\,.
\]
On the other hand, keeping in mind that $\Omega$ is a minimiser of the energy and by~(\ref{optsubdiv}) we obtain
\[\begin{split}
\F_\eps(E)&\geq \widetilde\F_\eps(E)
\geq \widetilde\F_\eps(\Omega)
=\sumn_{i=1}^H \F_\eps(\Omega^i)
=\F_\eps(\Omega^1) +\sumn_{i=2}^H \F_\eps(\Omega^i)
\geq \F_\eps(\Omega^1) +\sumn_{i=2}^H P(\Omega^i)\\
&\geq \F_\eps(\Omega^1) +N\omega_N^{1/N}\sumn_{i=2}^H |\Omega^i|^{\frac{N-1}N}
\geq \F_\eps(\Omega^1) +N\omega_N^{1/N}\, m^{\frac{N-1}N}\,.
\end{split}\]
Putting together the last two estimates, and recalling that $\F_\eps(\Omega^1)\leq \F_\eps(\Omega)\leq \F_\eps(B)$, we get
\[
N\omega_N^{1/N}\, m^{\frac{N-1}N} \leq \frac{3\F_\eps(B)}{\omega_N}\, m\,,
\]
which is impossible if $0<m<\overline m$ as soon as $\overline m$ is small enough. We deduce that necessarily $m=0$, so that actually $E=\Omega$ and we have obtained $\F_\eps(\Omega)=\widetilde\F_\eps(\Omega)$. In particular, $\Omega$ is also a minimiser of $\F_\eps$ and $\min\F_\eps=\min\widetilde\F_\eps$, and we only have to get~(\ref{connectedness}).

Let us then assume that~(\ref{connectedness}) is false, and let $E\subseteq\Omega$ be such that $P(\Omega)=P(E)+P(\Omega\setminus E)$, and $0<|E|<\Omega$. Then, since $g$ is strictly positive and thus
\[
\int_E \int_{\Omega\setminus E} g(y-x)\,dy\,dx>0\,,
\]
we have
\[
\widetilde\F_\eps(\Omega) \leq \F_\eps(E)+\F_\eps(\Omega\setminus E)
= P(E) + P(\Omega\setminus E) + \eps\Big(\En(E) +\En(\Omega\setminus E)\Big)
< P(\Omega) + \eps\En(\Omega) = \F_\eps(\Omega)\,,
\]
against the equality $\F_\eps(\Omega)=\widetilde\F_\eps(\Omega)$, that has been already established. This concludes the proof.
\end{proof}

Let us now pass to show that solution have no holes, for small $\eps$. In particular, they are simply connected if $N=2$.

\begin{lemma}[Solutions have no holes for $\eps\ll 1$]\label{noholes}
If $g$ is decreasing, then there exists $0<\eps_2<\eps_1$ such that, for every $0<\eps<\eps_2$, every minimiser of $\F_\eps$ (or, equivalently, of $\widetilde\F_\eps$, by Lemma~\ref{solconn}) has no holes. Equivalently, there is no set $G\subseteq \R^N\setminus \Omega$ with $|G|>0$ such that
\begin{equation}\label{imposs}
P(\Omega) = P(\Omega\cup G) + P(G)\,.
\end{equation}
\end{lemma}
\begin{proof}
Let us assume that $g$ is decreasing, that $\Omega\in\A$ is a minimiser of $\F_\eps$ for some $\eps<\eps_1$, and that there exists a set $G\subseteq\R^N\setminus\Omega$ with $m:=|G|>0$ and such that~(\ref{imposs}) holds. We have to find a contradiction if $\eps$ is smaller than a suitable $\eps_2$. First of all, we can show that $m\ll 1$. In fact, by the optimality of $\Omega$ we find
\[\begin{split}
N\omega_N + \eps\En(B)&= \F_\eps(B) \geq \F_\eps(\Omega) \geq P(\Omega)
= P(\Omega\cup G)+P(G)\\
&\geq N\omega_N^{1/N} \Big( (\omega_N+m)^{\frac{N-1}N} + m^{\frac{N-1}N}\Big)
\geq N\omega_N + N\omega_N^{1/N} m^{\frac{N-1}N}\,,
\end{split}\]
which gives
\begin{equation}\label{msmall}
m \leq \bigg(\frac{\En(B)}{N\omega_N^{1/N}}\bigg)^{\frac N{N-1}}\, \eps^{\frac N{N-1}}\,.
\end{equation}
Since $\eps< \eps_1$, in particular $m$ is bounded, so by Lemma~\ref{FvG} there exists a constant $C$, only depending on $g$ and on $\eps_1$ (thus ultimately only on $g$) such that
\[
\En(\Omega\cup G) - \En(\Omega) = \En(G) + 2 \En(\Omega,G) \leq C m\,.
\]
As a consequence, again by~(\ref{imposs}) we obtain
\begin{equation}\label{Feogsm}\begin{split}
\F_\eps(\Omega\cup G) &= P(\Omega\cup G) +\eps\En(\Omega\cup G)
= P(\Omega) - P(G) +\eps\En(\Omega\cup G)\\
&\leq P(\Omega) - N\omega_N^{1/N} m^{\frac{N-1}N} +\eps\En(\Omega) + \eps C m
=\F_\eps(\Omega) - N\omega_N^{1/N} m^{\frac{N-1}N}+\eps C m\\
&<\F_\eps(\Omega)\,,
\end{split}\end{equation}
where the last inequality holds as soon as $m$ is small enough, hence by~(\ref{msmall}) as soon as $\eps<\eps_2$ for a suitably small $\eps_2$.\par

In order to find a contradiction with the minimality of $\Omega$, we let $t\in\R$ be a number such that
\begin{align*}
|F|=\omega_N\,, && \hbox{where} && F = \Big\{x\in \Omega\cup G,\, x_1<t\Big\}\,.
\end{align*}
Notice that $\partial^* F\setminus \partial^*(\Omega\cup G)$ is contained in the hyperplane $t\times\R^{N-1}$. Moreover, for every $x'\in\R^{N-1}$ such that $(t,x')\in \partial^* F$ there exists some $x_1\geq t$ such that $(x_1,x')\in\partial^* (\Omega\cup G)$. Thus $P(F)\leq P(\Omega\cup G)$. Moreover, we also have $\En(F)\leq \En(\Omega\cup G)$ since $F\subseteq\Omega\cup G$. Hence, also by~(\ref{Feogsm}) we have
\[
\F_\eps(F) \leq \F_\eps(\Omega\cup G) < \F_\eps(\Omega)\,,
\]
which is the desired contradiction since $\Omega$ is optimal and $F\in\A$.
\end{proof}

We can now show that the Fraenkel asymmetry of optimal sets converge to $0$ when $\eps\searrow 0$, that is, the optimal sets converge in the $L^1$ sense to the unit ball, up to translations.

\begin{lemma}[Vanishing Fraenkel asymmetry]\label{lemmafraenkel}
For every $\delta>0$ there exists $\eps(\delta)>0$ such that, if $\Omega$ is a minimiser of $\F_\eps$ with some $\eps<\eps(\delta)$, then
\[
\min \Big\{ \big|\Omega\Delta (z+B)\big|,\, z\in \R^N\Big\} \leq \delta\,.
\]
\end{lemma}

\begin{proof}
By the quantitative isoperimetric inequality (see~\cite{FMP,CL}) there exists a geometric constant $C=C(N)$ such that, for every $\Omega\subseteq\R^N$ with $|\Omega|=\omega_N$, up to a translation we have
\[
P(\Omega)\geq P(B) + C|\Omega\Delta B|^2\,.
\]
If $\Omega$ is a minimiser of $\F_\eps$, we have then
\[
P(B) + \eps \En(B) = \F_\eps(B) \geq \F_\eps(\Omega) \geq P(\Omega) \geq P(B) +C|\Omega\Delta B|^2\,,
\]
from which we deduce
\[
|\Omega\Delta B| \leq \sqrt{\frac\eps C\, \En(B)}\,.
\]
The thesis then follows.
\end{proof}

The following observation will be useful in the proof of Theorem~\ref{2dcase}.

\begin{lemma}\label{lemadm}
A radial and decreasing function $g:\R^N\setminus\{0\} \to (0,+\infty)$ is admissible, in the sense of Definition~\ref{def:admissible}, if and only if
\begin{equation}\label{N-1}
\int_0^1 g(t)t^{N-1}\,dt<+\infty\,.
\end{equation}
Under this assumption, the function $\Phi:[0,+\infty)\to [0,+\infty]$ defined as
\begin{equation}\label{defPhi}
\Phi(t) := \int_B g(y-x)\, dx\qquad
\text{for }y\in\partial B_t\,,
\end{equation}
is decreasing and locally Lipschitz continuous in $\R^+\setminus \{1\}$. Moreover, $\Phi$ is locally Lipschitz continuous on the whole $\R^+$ if and only if
\begin{equation}\label{N-2}
\int_0^1 g(t)t^{N-2}\,dt<+\infty\,.
\end{equation}
\end{lemma}

\begin{proof}
The first property simply follows by integration in polar coordinates. Indeed, assume that $g$ is radial and decreasing. Then, on one hand one has
\[\begin{split}
\En(B)&=\int_B \int_B g(y-x)\,dy\,dx
=\int_B \int_{-x+B} g(y)\,dy\,dx
\geq \int_{B_{1/2}} \int_{B_{1/2}} g(y)\,dy\,dx\\
&= \frac{\omega_N}{2^N}\, \int_{t=0}^{1/2} t^{N-1}\,dt\,,
\end{split}\]
and on the other hand
\[
\En(B)=\int_B \int_{-x+B} g(y)\,dy\,dx\leq \int_B \int_{B_2} g(y)\,dy\,dx
=\omega_N \int_{t=0}^2 t^{N-1}\,dt\,.
\]
Since $g$ is decreasing, $\int_0^\tau g(t)\,dt$ is either finite for every $\tau>0$ or infinite for every $\tau>0$, thus the first property follows.\par
Let us then assume that the property~(\ref{N-1}) is satisfied, and let $\Phi$ be the function defined in~(\ref{defPhi}). The fact that $\Phi$ is decreasing and smooth in $\R^+\setminus\{1\}$, so in particular locally Lipschitz in $\R^+\setminus\{1\}$, directly follows from the definition. Concerning the Lipschitz property around $1$, it clearly holds if and only if there exists a constant $C>0$ such that
\[
\int_{R_\eps} g(|x|)dx<C\eps\, ,
\]
where $R_\eps = \big\{x=(x',x_N)\in\R\times\R^{N-1}:\, |x'|<1,\, |x_N|<\eps\big\}$. A simple integration in polar coordinates ensures that this is equivalent to~(\ref{N-2}).
\end{proof}

\subsection{The $2$-dimensional case: minimality of the ball for small $\eps$}

We now restrict ourselves to the $2$-dimensional case $N=2$. Notice that, in this case, the inequality~(\ref{N-2}) reduces to~(\ref{ipLip}). For further use, we show the following geometric estimate.
\begin{lemma}\label{lempag1}
Let $0<\bar\theta<\pi/2$ and let $-\cos(\bar\theta)/8<\delta<\cos(\bar\theta)/8$. As in Figure~\ref{Fignotat}, let $\tau(0)$ be the length of the arc of circle of radius $1$, centred at $O\equiv (0,0)$ and connecting the points $P\equiv (\cos\bar\theta,\sin\bar\theta)$ and $Q\equiv (\cos\bar\theta,-\sin\bar\theta)$ through $(1,0)$, and let $\tau=\tau(\delta)$ be the length of the arc of circle connecting the points $P$ and $Q$ and passing through $S\equiv (1+\delta,0)$. Then,
\begin{equation}\label{estpag1}
\tau(\delta)- \tau(0) \geq \mu + \frac{\cos \bar\theta}6\, \delta\mu\,,
\end{equation}
where $\mu=\mu(\delta)$ is the signed area enclosed between the two arcs of circle (positive for $\delta>0$).
\end{lemma}
\begin{proof}
Let $0<\bar\theta<\pi/2$ be fixed. For every $-\cos(\bar\theta)/8<\delta<\cos(\bar\theta)/8$, let us call $R$ the centre of the arc $\tau(\delta)$, and let us call $\eta=\eta(\delta)$ the number such that $R\equiv (\eta,0)$. Let moreover $\rho=\rho(\delta)$ be the radius of the arc $\tau(\delta)$, hence the length of the segment $PR$, and $\theta=\theta(\delta)$ the angle $\angle PRS$, so that $\theta(0)=\bar\theta$. The different quantities are depicted in Figure~\ref{Fignotat}.\par
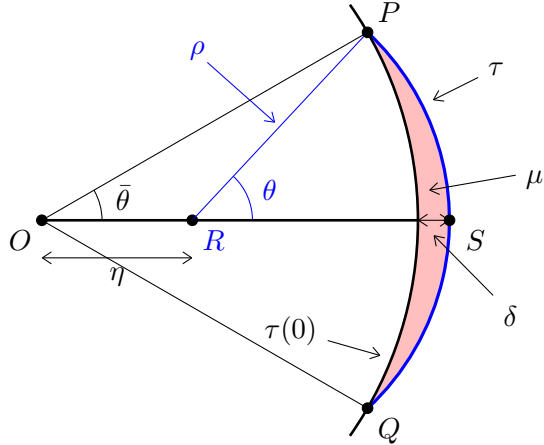
\begin{figure}[htbp]
\begin{tikzpicture}[>=>>>]
\filldraw[fill=red!25!white, draw=black, line width=0] (4.33,2.5) arc(47.02:-47.02:3.42) arc(-30:30:5);
\filldraw (2,0) circle(2pt);
\draw[blue] (2,0) node[anchor=north west] {$R$};
\draw (.8,0) arc(0:30:.8);
\draw (.85,.3) node[anchor=west] {$\bar\theta$};
\draw [blue] (2,0) -- (4.33,2.5);
\draw [blue,->] (2.2,2.1) -- (3.1,1.3);
\draw[blue] (2.3,2) node[anchor=south east] {$\rho$};
\draw[blue] (2.8,0) arc(0:47.02:.8);
\draw[blue] (2.8,.4) node[anchor=west] {$\theta$};
\draw[blue, line width=1.2] (4.33,2.5) arc(47.02:-47.02:3.42);
\filldraw (0,0) circle(2pt);
\draw (0,0) node[anchor=north east] {$O$};
\draw[line width=1] (0,0) -- (5,0) arc(0:35:5);
\draw[line width=1] (5,0) arc(0:-35:5);
\filldraw (4.33,2.5) circle(2pt);
\filldraw (4.33,-2.5) circle(2pt);
\draw (4.33,2.5) node[anchor=south west] {$P$};
\draw (4.33,-2.5) node[anchor=north west] {$Q$};
\draw (4.33,2.5) -- (0,0) -- (4.33,-2.5);
\draw [<->] (0,-.5) -- (2,-.5);
\draw (1,-.5) node[anchor=north] {$\eta$};
\draw (5.8,1.8) node[anchor=south west] {$\tau$};
\draw [->] (5.8,1.9) -- (5.2,1.6);
\draw (6.3,.8) node[anchor=north west] {$\mu$};
\draw [->] (6.3,.6) -- (5.15,.4);
\draw (3.8,-1.8) node[anchor=south east] {$\tau(0)$};
\draw [->] (3.8,-1.6) -- (4.5,-1.8);
\filldraw (5.42,0) circle(2pt);
\draw (5.5,0) node[anchor=north west] {$S$};
\draw [<->] (5,0) -- (5.42,0);
\draw [->] (6,-1) -- (5.21,-.2);
\draw (6,-1) node[anchor=north west] {$\delta$};
\end{tikzpicture}
\caption{Situation for the proof of Lemma~\ref{lempag1}.}\label{Fignotat}
\end{figure}
A few elementary trigonometric calculations give that
\[
\rho=\sqrt{(\cos\bar\theta-\eta)^2+\sin^2\bar\theta} = 1 -\eta\cos\bar\theta+o(\eta)\,,
\]
and then, since $\sin\bar\theta=\rho\sin\theta$, that
\[
\theta=\bar\theta+\eta\sin\bar\theta + o(\eta)\,.
\]
Keeping in mind that $1+\delta=\eta+\rho$, we deduce that
\[
\eta = \frac 1{1-\cos\bar\theta}\,\delta + o(\delta)\,,
\]
thus the above first order expansions ensure that
\begin{align}\label{rho'theta'0}
\rho'(0) = -\frac{\cos\bar\theta}{1-\cos\bar\theta}\,, && \theta'(0) = \frac{\sin\bar\theta}{1-\cos\bar\theta}\,.
\end{align}
Concerning the area $\mu$, one has by construction
\[
\mu=\rho^2\theta+\eta\sin\bar\theta-\bar\theta=
2\, \frac{\sin\bar\theta-\bar\theta\cos\bar\theta}{1-\cos\bar\theta}\, \delta + o(\delta)\,.
\]
Notice that, of course, $\mu$ is positive if and only if so is $\delta$. Finally, it is easy to evaluate $\tau$ as
\[
\tau= 2\rho\theta = 2\bar\theta +2 \,\frac{\sin\bar\theta-\bar\theta\cos\bar\theta}{1-\cos\bar\theta}\, \delta + o(\delta)
=\tau(0) +2 \,\frac{\sin\bar\theta-\bar\theta\cos\bar\theta}{1-\cos\bar\theta}\, \delta + o(\delta)\,,
\]
so the last two estimates give
\begin{equation}\label{mu'tau'0}
\mu'(0)= 2\, \frac{\sin\bar\theta-\bar\theta\cos\bar\theta}{1-\cos\bar\theta}= \tau'(0)\,.
\end{equation}
It is now easy to evaluate the derivatives of the functions $\rho,\, \theta,\, \mu$ and $\tau$ at values of $\delta$ different from $0$. Indeed, since of course the lengths $\tau,\, \delta$ and $\rho$ are linear with respect to the radius, the area $\mu$ is quadratic, and the angle $\theta$ is $0$-homogeneous, from~(\ref{rho'theta'0}) and~(\ref{mu'tau'0}) we directly have
\begin{align}\label{derivatives}
\rho'= -\frac{\cos\theta}{1-\cos\theta}\,, && \theta' = \frac{\sin\theta}{\rho(1-\cos\theta)}\,, &&
\tau'=2\, \frac{\sin\theta-\theta\cos\theta}{1-\cos\theta}\,, && \mu'=\rho\tau'\,.
\end{align}
In order to show~(\ref{estpag1}), we will now argue separately for the case $\delta>0$ and $\delta<0$.

\case{I}{The case $\delta>0$.}
Let us start considering the case when $\delta>0$. We claim that
\begin{equation}\label{claimcase1}
\rho(\sigma)\leq 1-\frac {\cos\bar\theta}2\,\sigma\qquad \forall\ 0<\sigma<\frac{\cos\bar\theta}8\,.
\end{equation}
To show this estimate, we call $\delta_1>0$ the value of $\delta$ such that the corresponding $\eta$ is $\eta(\delta_1)=(\cos\bar\theta)/2$. Notice that, as a consequence, one has
\begin{equation}\label{hererhodelta1}
\rho(\delta_1)=\sqrt{\frac{\cos^2\bar\theta}4+\sin^2\bar\theta} = \frac 12\, \sqrt{1+3\sin^2\bar\theta} \geq \frac{1+\sin^2\bar\theta}2\,,
\end{equation}
and since $(\cos\bar\theta)/2+\rho(\delta_1)=1+\delta_1$ we get
\begin{equation}\label{estdel1}
\delta_1 = \frac{\cos\bar\theta}2 + \rho(\delta_1) -1 \geq \frac{\cos\bar\theta+\sin^2\bar\theta-1}2
=\frac{\cos\bar\theta(1-\cos\bar\theta)}2\,.
\end{equation}
Notice now that, for every $0<\sigma<\delta_1$, by~(\ref{derivatives}) one has
\[
\rho'(\sigma) = -\frac{\cos(\theta(\sigma))}{1-\cos(\theta(\sigma))} \leq -\cos(\theta(\sigma)) \leq -\rho(\sigma)\cos(\theta(\sigma)) \leq -\frac{\cos\bar\theta}2\,,
\]
from which the inequality $\rho(\sigma)\leq 1 - \sigma\cos\bar\theta/2$ follows for every $0<\sigma\leq \delta_1$. We have then already obtained~(\ref{claimcase1}) if $\delta_1\geq\cos\bar\theta/8$, which by~(\ref{estdel1}) is surely true if $\cos\bar\theta\leq 3/4$.\par

Let us instead assume that $\cos\bar\theta> 3/4$. In this case notice that, by elementary geometric reasons, $\rho(\sigma)\leq \rho(\delta_1)$ for every $\delta_1\leq \sigma\leq \delta_2$, where $\delta_2$ is such that the corresponding $\eta$ is $\eta(\delta_2)=3(\cos\bar\theta)/2$. Since
\[
\delta_2=\eta(\delta_2)+\rho(\delta_2)-1\geq 3\cos\bar\theta-1\geq \frac{\cos\bar\theta}8\,,
\]
we deduce~(\ref{claimcase1}) also with $\cos\bar\theta>3/4$, since the inequality has already been proved for $0<\sigma\leq \delta_1$, and for every $\delta_1<\sigma<\cos\bar\theta/8$ one has, also by~(\ref{hererhodelta1}),
\[
\rho(\sigma)\leq \rho(\delta_1) = \frac 12\, \sqrt{1+3\sin^2\bar\theta} \leq 1 - \frac{\cos\bar\theta}4 \leq 1 - \frac{\cos\bar\theta}2\, \sigma\,.
\]
Making use of~(\ref{claimcase1}) we can easily obtain~(\ref{estpag1}). In fact, since $0<\delta<\cos\bar\theta/8$, by~(\ref{derivatives}) and~(\ref{claimcase1}) we have, for $0<\sigma<\delta$,
\[
\mu'(\sigma) = \rho(\sigma) \tau'(\sigma) \leq \tau'(\sigma) - \frac{\cos\bar\theta}2\,\sigma \tau'(\sigma)\,.
\]
Moreover, a simple calculation ensures that the function
\[
\theta\mapsto \frac{\sin\theta-\theta\cos\theta}{1-\cos\theta}
\]
is increasing, hence by~(\ref{derivatives}) also $\tau'$ is an increasing function of $\delta$. Therefore, we can estimate
\[\begin{split}
\tau(\delta)-\tau(0)&= \int_{\sigma=0}^\delta \tau'(\sigma)\,d\sigma
\geq \int_0^\delta \mu'(\sigma)\,d\sigma + \int_0^\delta \frac{\cos\bar\theta}2 \sigma\tau'(\sigma)\,d\sigma\\
&= \mu + \frac{\cos\bar\theta}2\int_0^\delta \sigma\tau'(\sigma)\,d\sigma
\geq \mu + \frac{\cos\bar\theta}4\, \delta \int_0^\delta \tau'(\sigma)\,d\sigma\\
&\geq \mu + \frac{\cos\bar\theta}4\, \delta \int_0^\delta \rho(\sigma)\tau'(\sigma)\,d\sigma
=\mu + \frac{\cos\bar\theta}4\, \delta \int_0^\delta \mu'(\sigma)\,d\sigma
=\mu + \frac{\cos\bar\theta}4\, \delta \mu\,,
\end{split}\]
which is stronger than~(\ref{estpag1}).

\case{II}{The case $\delta<0$.}
Let us now consider the case when $\delta<0$. In this case, we call $\delta_1=\cos\bar\theta-1$. Notice that $\delta_1$ is negative, and it corresponds to the situation in which $S$ is the middle point of the segment $PQ$, hence in particular $\rho(\sigma)\to +\infty$ for $\sigma\searrow \delta_1$. We first aim to show that
\begin{equation}\label{mainestcase2}
\mu''(\sigma)\geq 0\qquad \forall \,\delta_1<\sigma<0\,.
\end{equation}
In fact, by homogeneity, it is enough to show that $\mu''(0)\geq 0$. Keeping in mind~(\ref{derivatives}), we have
\[
\mu''(0) = (\rho\tau')'(0) = 2\bigg(-\frac{\cos\bar\theta(\sin\bar\theta-\bar\theta\cos\bar\theta)}{(1-\cos\bar\theta)^2} + \frac{\bar\theta\sin^2\bar\theta-\sin^3\bar\theta}{(1-\cos\bar\theta)^3}\bigg)\,,
\]
so that we are reduced to check that
\[
\bar\theta\sin^2\bar\theta-\sin^3\bar\theta \geq \cos\bar\theta(\sin\bar\theta-\bar\theta\cos\bar\theta)(1-\cos\bar\theta)\,,
\]
which in turn can be rewritten as
\[
\bar\theta(1+\cos\bar\theta+\cos^2\bar\theta) \geq \sin\bar\theta(1+2\cos\bar\theta)\,.
\]
In other words, we have to show $f(\bar\theta)\geq 0$, where $f(x)=x(1+\cos x + \cos^2 x) -\sin x(1+2\cos x)$. Since $0<\bar\theta<\pi/2$ and $f(0)=0$, it is sufficient to show that $f'(x)\geq 0$ for every $0<x<\pi/2$, and this is equivalent to say that
\[
3\sin x \geq x(1+2\cos x) \qquad \forall\, 0<x<\frac\pi 2\,.
\]
Once again, we can call $\tilde f(x) =3\sin x - x(1+2\cos x)$, observe that $\tilde f(0)=0$, and then it is enough to show that $\tilde f'(x)\geq 0$ for every $0<x<\pi/2$. And finally,
\[
\tilde f'(x) = \cos x -1 +2x\sin x\geq \cos^2 x-1 +2x\sin x = \sin x (2x - \sin x) \geq 0\qquad \forall\, 0<x<\frac \pi 2\,,
\]
so~(\ref{mainestcase2}) is proved.\par

A simple trigonometric calculation ensures that
\begin{equation}\label{boundmu'd}
\mu'(\delta_1)=\frac 43\, \sin\bar\theta\,,
\end{equation}
hence for every $\delta_1\leq \sigma\leq 0$ we have
\begin{equation}\label{boundmu'}
\frac 23\, \mu'(0) \leq \mu'(\sigma)\leq \mu'(0)\,,
\end{equation}
where the second inequality is a direct consequence of~(\ref{mainestcase2}), and the first follows by~(\ref{mainestcase2}) since
\[
\mu'(\sigma)\geq \mu'(\delta_1) = \frac 43\, \sin\bar\theta \geq
\frac43\, \frac{\sin\bar\theta-\bar\theta\cos\bar\theta}{1-\cos\bar\theta} =\frac 23\, \mu'(0)\,.
\]
We can then argue more or less as in Case~I. In fact, observe that for every $\delta_1\leq \sigma\leq 0$ one has $0\leq \theta(\sigma)\leq \bar\theta$, hence by~(\ref{derivatives})
\[
\rho'(\sigma)=-\frac{\cos(\theta(\sigma))}{1-\cos(\theta(\sigma))}\leq -\frac{\cos\bar\theta}{1-\cos\bar\theta}\leq -\cos\bar\theta\,,
\]
so that
\[
\frac 1{\rho(\sigma)} \leq \frac 1{1+\cos\bar\theta|\sigma|} \leq 1 - \frac{\cos\bar\theta}2\,|\sigma|\,.
\]
Hence, by~(\ref{derivatives}) and~(\ref{boundmu'}), for every $\delta_1\leq \delta\leq 0$ we have
\[\begin{split}
\tau(\delta)-\tau(0) &= - \int_{\sigma=\delta}^0 \tau'(\sigma)\,d\sigma
= - \int_\delta^0 \frac{\mu'(\sigma)}{\rho(\sigma)}\,d\sigma
\geq \mu +\frac{\cos\bar\theta}2\, \int_\delta^0 \mu'(\sigma) |\sigma|\,d\sigma\\
&\geq \mu +\frac{\cos\bar\theta}3\, \int_\delta^0 \mu'(0) |\sigma|\,d\sigma
=\mu -\frac{\cos\bar\theta}6\,\delta \int_\delta^0 \mu'(0) \,d\sigma
\geq \mu -\frac{\cos\bar\theta}6\,\delta \int_\delta^0 \mu'(\sigma) \,d\sigma\\
&= \mu +\frac{\cos\bar\theta}6\, \delta \mu\,,
\end{split}\]
which is~(\ref{estpag1}).\par
We have then proved~(\ref{estpag1}) for every $\delta_1\leq \delta\leq 0$, while we have to prove it for every $-\cos\bar\theta/8< \delta< 0$. However, the missing cases are very easy. In fact, if $\delta\leq 2\delta_1$, then an immediate geometric argument ensures that $\tau(\delta)\geq \tau(0)$. Hence, for every $-\cos\bar\theta/8<\delta\leq 2\delta_1$ we clearly have~(\ref{estpag1}) since
\[
\tau(\delta)-\tau(0)\geq 0 \geq \mu + \frac{\cos\bar\theta}6\,\delta\mu\,,
\]
where the last inequality is true since $\mu\leq 0$ and
\[
1+\frac{\cos\bar\theta}2\,\delta \geq 1-\frac{\cos^2\bar\theta}{16}>0\,. 
\]
And finally, to show~(\ref{estpag1}) for $2\delta_1<\delta<\delta_1$, it is enough to check that $\Psi'(\sigma)>0$ for every $2\delta_1<\sigma<0$, where
\[
\Psi(\sigma) = \mu(\sigma) + \frac{\cos\bar\theta}6\,\sigma\mu(\sigma)\,.
\]
Indeed, once we have proved $\Psi'>0$ in the interval $[2\delta_1,\delta_1]$, (\ref{estpag1}) immediately follows since for every $2\delta_1<\delta<\delta_1$
\[
\tau(\delta)-\tau(0) \geq \tau(\delta_1)-\tau(0) \geq \Psi(\delta_1) \geq \Psi(\delta) = \mu + \frac{\cos\bar\theta}6\,\sigma\mu\,.
\]
Let us then show the inequality $\Psi'>0$. Since for every $2\delta_1<\sigma<0$ we have by trigonometric reasons and by~(\ref{boundmu'}) and~(\ref{boundmu'd}) that
\begin{align*}
-2\big(\bar\theta-\sin\bar\theta\cos\bar\theta\big)\leq \mu(\sigma)\leq -\big(\bar\theta-\sin\bar\theta\cos\bar\theta\big)\,, &&
\mu'(\sigma)\geq \frac 43\,\sin\bar\theta\,,
\end{align*}
keeping in mind~(\ref{derivatives}) we have
\[
\Psi'(\sigma) = \mu'(\sigma)\bigg(1+\frac{\cos\bar\theta}6\,\sigma\bigg) + \frac{\cos\bar\theta}6\,\mu(\sigma)
\geq \frac 23\,\sin\bar\theta - \frac{\bar\theta-\sin\bar\theta\cos\bar\theta}3
\geq \frac {2\sin\bar\theta-\bar\theta}3>0\,,
\]
so as observed before the proof is completed.
\end{proof}
Let us now consider a set $E\in\A$. Up to a translation, we can assume that
\begin{equation}\label{centered}
\min \Big\{ \big|E\Delta (z+B)\big|,\, z\in \R^N\Big\} = |E\Delta B|\,,
\end{equation}
since the minimum exists by a simple compactness argument. For any set $E\in\A$, under the assumption~(\ref{centered}), we call
\begin{equation}\label{defdeltaA}
\nu:=\frac{|E\Delta B|}2=|E\setminus B|=|B\setminus E|\,,
\end{equation}
and
\begin{align}\label{precdelta}
\delta^+:= \sup \big\{ s>0:\, |E\setminus B_s|>0\big\} - 1 \,,&&
\delta^-:=1 - \inf\big\{s:\, |B_s\setminus E|>0\big\}\,.
\end{align}
We want now to show the following result.
\begin{lemma}\label{AdelC}
Let $E\in\A$ be a set satisfying~(\ref{centered}), connected and with no holes in the sense of~(\ref{connectedness}) and~(\ref{imposs}). Let moreover $\delta^\pm$ and $\nu$ be defined as in~(\ref{defdeltaA}) and~(\ref{precdelta}). Then,
\begin{equation}\label{thesisquad}
P(E)-P(B) \geq \frac 1C\,\nu(\delta^++\delta^-)\,,
\end{equation}
for some purely geometric constant $C$.
\end{lemma}
\begin{proof}
Since the proof is quite involved, we divide it in some steps. Moreover, we will consider separately the situation in $E\setminus B$ and in $B\setminus E$.
\step{I}{The families $\Gamma,\, \Gamma^{t,\beta},\, \Gamma_{t,\beta}$}
We consider a class of ``generalised possible boundaries'' as follows. We let $\Gamma_0$ be the class of all ${\rm C}^1$, injective curves $\gamma:\S^1\to\R^2$, parametrized with constant speed and counterclockwise (that is, $|\gamma'|$ is constant and all the points internal to the curve have degree $1$ with respect to the curve itself). We call then $\Gamma$ the class of functions $\gamma:\S^1\to\R^2$ which are uniform limits of elements of $\Gamma_0$. Notice that any curve in $\Gamma$ is Lipschitz continuous and with $|\gamma'|$ constant, but it is not necessarily injective, hence it is not necessarily a Jordan curve. Given a curve $\gamma\in\Gamma$, and writing for brevity $\gamma$, with a small abuse of notation, also to denote the set $\gamma(\S^1)\subseteq\R^2$, every point $x\in\R^2\setminus \gamma$ has either degree $0$ or degree $1$, and the set $E_\gamma\subseteq\R^2$ of points with degree $1$ is a bounded, open set. Notice that a same set $E_\gamma$ corresponds to different curves $\gamma$, even up to rotations of $\S^1$. Indeed, the set $E_\gamma$ does not change if one adds to a curve $\gamma$ a Lipschitz curve in $\R^2\setminus E_\gamma$, with one endpoint in $\gamma$, percurred once outwards and then once inwards (see for instance the fourth picture in Figure~\ref{FigIsop}). For any positive constants $t$ and $\beta$, we define
\begin{align*}
\Gamma^{t,\beta} := \Big\{ \gamma\in\Gamma:\, E_\gamma\supseteq B,\, \max\big\{|\gamma(t)|,\, t\in\S^1\big\}=1+t\,, |E_\gamma\setminus B|=\beta\Big\}\,,\\
\Gamma_{t,\beta} := \Big\{ \gamma\in\Gamma:\, E_\gamma\subseteq B,\, \min\big\{|\gamma(t)|,\, t\in\S^1\big\}=1-t\,, |B\setminus E_\gamma|=\beta\Big\}\,.
\end{align*}
Notice that the family $\Gamma^{t,\beta}$ (resp., $\Gamma_{t,\beta}$) is non-empty only if $\pi+\beta\leq \pi(1+t)^2$ (resp., $t\leq 1$ and $\pi-\beta\geq \pi(1-t)^2$). Finally, we define the ``length'' of every curve $\gamma\in\Gamma$ as
\[
\ell(\gamma) = \int_{\S^1} |\gamma'(t)|\,dt\,.
\]
\step{II}{The curves $\bar\gamma^{t,\beta}$ and $\bar\gamma_{t,\beta}$.}
Let $t$ and $\beta$ be positive numbers such that $\Gamma^{t,\beta}$ (resp., $\Gamma_{t,\beta}$) is non-empty. By a simple compactness argument, this family contains a curve $\bar\gamma^{t,\beta}$ (resp., $\bar\gamma_{t,\beta}$) with minimal length. It is quite standard to describe these curves, depending on $t$ and $\beta$.\par

Let us start with the curve $\bar\gamma^{t,\beta}$. The ``free boundary'' (i.e. the set of points $x\in \bar\gamma^{t,\beta}$ with $1<|x|<1+t$) has to be made by arcs of circle, all of the same radius. Moreover, these arcs have to meet $\partial B$ and each other tangentially, except at points $x$ with $|x|=1+t$. As a consequence, one readily derives that the situation is one of the four of Figure~\ref{FigIsop}, which are depicted for increasing values of $t$. More precisely, for each given $\beta\geq 0$, if $t$ is small enough (but still such that $\Gamma^{t,\beta}\neq \emptyset$) then $\bar\gamma^{t,\beta}\setminus\partial B$ is done by two arcs with positive curvature, meeting with a corner, as in the first picture in the left. If $t$ increases then the curvature of the two arcs decreases, and it becomes null for some $t$ (second picture, where the two arcs are actually segments) and then negative (third picture). Eventually, for $t$ large enough, the two arcs of the curve $\bar\gamma^{t,\beta}$ meet tangentially at some distance $t_{\rm ext}<t$ from $\partial B$, and then there is a final segment of length $t-t_{\rm ext}$, so to reach the final distance $t$. Notice that the final segment counts twice in $\ell(\bar\gamma^{t,\beta})$, and in particular this shows that both the inequalities $\ell(\gamma)\geq\H^1(\gamma)\geq P(E_\gamma)$ for elements $\gamma\in\Gamma$ can be strict. For ease of notation later, we set $t_{\rm ext}=t$ in the first three cases.

\begin{figure}[htbp]
\begin{tikzpicture}
\filldraw[fill=red!25!white, draw=black, line width=0] (1.03,1.29) arc(40:18.75:4) arc(-18.75:-40:4) arc(-40:40:2);
\filldraw (-0.5,0) circle(1.5pt);
\draw (-0.5,0) node[anchor=north east] {$O$};
\draw[dashed] (-.5,0) -- (1.76,0);
\draw[blue, line width=1] (1.03,1.29) arc (40:18.75:4);
\draw[blue, line width=1] (1.03,-1.29) arc (-40:-18.75:4);
\draw[line width=1] (1.5,0) arc (0:45:2);
\draw[line width=1] (1.5,0) arc (0:-45:2);
\filldraw (1.76,0) circle (1.5pt);
\draw (-0.5,0) -- (1.03,1.29);
\draw (0,0) arc(0:40:.5);
\draw (0,0.25) node[anchor=west] {$\theta$};
\draw (1.76,0) node[anchor=south west] {$\bar x$};
\draw[<->] (1.5,-1) -- (1.76,-1);
\draw (1.63,-1) node[anchor=north] {$t$};
\draw (1.76,.8) node[anchor=south west] {$\beta$};
\draw [->>>] (1.86,.9) -- (1.55,.2);
\filldraw[fill=red!25!white, draw=black, line width=0] (5.23,1) -- (5.81,0) -- (5.23,-1) arc(-30:30:2);
\filldraw (3.5,0) circle(1.5pt);
\draw[blue, line width=1] (5.23,1) -- (5.81,0) -- (5.23,-1);
\draw (3.5,0) node[anchor=north east] {$O$};
\draw[line width=1] (5.5,0) arc (0:45:2);
\draw[line width=1] (5.5,0) arc (0:-45:2);
\draw[dashed] (3.5,0) -- (5.81,0);
\filldraw (5.79,0) circle (1.5pt);
\draw (3.5,0) -- (5.23,1);
\draw (4,0) arc(0:30:.5);
\draw (4,0.25) node[anchor=west] {$\theta$};
\draw (5.79,0) node[anchor=south west] {$\bar x$};
\draw[<->] (5.5,-1) -- (5.79,-1);
\draw (5.64,-1) node[anchor=north] {$t$};
\draw (5.76,.8) node[anchor=south west] {$\beta$};
\draw [->>>] (5.86,.9) -- (5.55,.2);
\filldraw[fill=red!25!white, draw=black, line width=1] (9.38,.68) arc(200:223.16:2) arc(-223.16:-200:2) arc (-20:20:2);
\filldraw (7.5,0) circle(1.5pt);
\draw (7.5,0) node[anchor=north east] {$O$};
\draw[dashed] (7.5,0) -- (9.8,0);
\draw[blue, line width=1] (9.38,.68) arc (200:223.16:2);
\draw[blue, line width=1] (9.38,-.68) arc (-200:-223.16:2);
\draw[line width=1] (9.5,0) arc (0:45:2);
\draw[line width=1] (9.5,0) arc (0:-45:2);
\filldraw (9.8,0) circle (1.5pt);
\draw (7.5,0) -- (9.38,.68);
\draw (8,0) arc(0:20:.5);
\draw (7.55,0.3) node[anchor=west] {$\theta$};
\draw (9.8,0) node[anchor=south west] {$\bar x$};
\draw[<->] (9.5,-1) -- (9.8,-1);
\draw (9.65,-1) node[anchor=north] {$t$};
\draw (9.76,.8) node[anchor=south west] {$\beta$};
\draw [->>>] (9.86,.9) -- (9.58,.1);
\filldraw[fill=red!25!white, draw=black, line width=1] (13.43,.52) arc(195:230.94:1) arc(-230.94:-195:1) arc (-15:15:2);
\filldraw (11.5,0) circle(1.5pt);
\draw (11.5,0) node[anchor=north east] {$O$};
\draw[dashed] (11.5,0) -- (13.77,0);
\draw[blue, line width=1] (13.43,.52) arc (195:230.94:1) -- (14,0);
\draw[blue, line width=1] (13.43,-.52) arc (-195:-230.94:1);
\draw[line width=1] (13.5,0) arc (0:45:2);
\draw[line width=1] (13.5,0) arc (0:-45:2);
\filldraw (14,0) circle (1.5pt);
\draw (11.5,0) -- (13.43,.52);
\draw (12,0) arc(0:15:.5);
\draw (11.55,0.3) node[anchor=west] {$\theta$};
\draw (14,0) node[anchor=south west] {$\bar x$};
\draw[<->] (13.5,-1) -- (14,-1);
\draw[<->] (13.5,1) -- (13.77,1);
\draw (13.75,-1) node[anchor=north] {$t$};
\draw (13.63,1.05) node[anchor=south] {$t_{\rm ext}$};
\draw (12.8,-.6) node[anchor=north east] {$\beta$};
\draw [->>>] (12.8,-.6) -- (13.4,-.1);
\end{tikzpicture}
\caption{The possible curves $\bar\gamma^{t,\beta}$ for increasing values of $t$.}\label{FigIsop}
\end{figure}
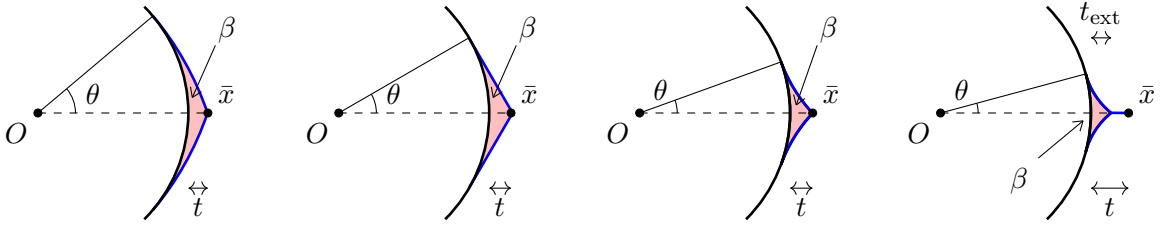
Let us now consider the curve $\bar\gamma_{t,\beta}$. The situation is similar to the preceding one, since again the free boundary, made by points $x\in \bar\gamma_{t,\beta}$ with $1-t<|x|<1$, is done by arcs of circle, all with the same radius, and meeting $\partial B$ and each other tangentially except at points $x$ with $|x|=1-t$. In this case the curvature of the free boundary has to be positive, so there are only two possibilities, depicted in Figure~\ref{Figdentro}. Namely, for $\beta$ fixed and $t$ such that $\Gamma_{t,\beta}\neq\emptyset$, if $t$ is smaller than some threshold then $\bar\gamma_{t,\beta}\setminus\partial B$ is done by two arcs, meeting with a corner at some point $\bar x$ with $|\bar x|=1-t$, as in the first picture. Instead, if $t$ is larger than the threshold, then $\bar\gamma_{t,\beta}\setminus\partial B$ is done by two arcs which meet tangentially at some point $\bar y$ at distance $t_{\rm int}$ from $\partial B$, plus a segment of length $t-t_{\rm int}$ so to reach distance $t$ from $\partial B$, as in the second picture. As before, for ease of notation later we set $t_{\rm int}=t$ in the first case.\par

\begin{figure}[htbp]
\begin{tikzpicture}
\filldraw[fill=red!25!white, draw=black, line width=0] (2.3,1.93) arc(40:-57.2:1.3) arc(57.2:-40:1.3) arc(-40:40:3);
\filldraw (0,0) circle(2pt);
\draw (0,0) node[anchor=north east] {$O$};
\draw (0,0) -- (2.3,1.93);
\draw (0.5,0) arc(0:40:.5);
\draw (0.5,0.25) node[anchor=west] {$\theta$};
\draw[blue, line width=1.2] (2.3,1.93) arc (40:-57.2:1.3);
\draw[blue, line width=1.2] (2.3,-1.93) arc (-40:57.2:1.3);
\draw[line width=1.2] (3,0) arc(0:45:3);
\draw[line width=1.2] (3,0) arc(0:-45:3);
\filldraw (2,0) circle(2pt);
\draw (2,0) node[anchor=north east] {$\bar x$};
\draw[dashed] (0,0) -- (2,0);
\draw[<<<->>>] (2,-1.2) -- (3,-1.2);
\draw (2.3,-1.2) node[anchor=north] {$t$};
\draw (3.5,1.1) node[anchor=west] {$\beta$};
\draw [->>>] (3.5,1) -- (2.7,0.1);
\filldraw[fill=red!25!white, draw=black, line width=0] (9.6,1.5) arc(30:-90:1) arc(90:-30:1) arc(-30:30:3);
\filldraw (7,0) circle(2pt);
\draw (7,0) node[anchor=north east] {$O$};
\draw (7,0) -- (9.6,1.5);
\draw (7.5,0) arc(0:30:.5);
\draw (7.5,0.25) node[anchor=west] {$\theta$};
\draw[blue, line width=1.2] (8.73,0) -- (8.2,0);
\draw[blue, line width=1.2] (9.6,1.5) arc (30:-90:1);
\draw[blue, line width=1.2] (9.6,-1.5) arc (-30:90:1);
\draw[line width=1.2] (10,0) arc(0:45:3);
\draw[line width=1.2] (10,0) arc(0:-45:3);
\filldraw (8.2,0) circle(2pt);
\draw (8.2,0) node[anchor=north east] {$\bar x$};
\filldraw (8.7,0) circle(2pt);
\draw (8.7,0) node[anchor=north west] {$\bar y$};
\draw[dashed] (7,0) -- (8.2,0);
\draw[<<<->>>] (8.2,-1.2) -- (10,-1.2);
\draw (9.1,-1.2) node[anchor=north] {$t$};
\draw[<<<->>>] (8.7,1.2) -- (10,1.2);
\draw (9.3,1.2) node[anchor=north] {$t_{\rm int}$};
\draw (10.5,1.1) node[anchor=west] {$\beta$};
\draw [->>>] (10.5,1) -- (9.7,0.1);
\end{tikzpicture}
\caption{The possible curves $\bar\gamma_{t,\beta}$ for increasing values of $t$.}\label{Figdentro}
\end{figure}
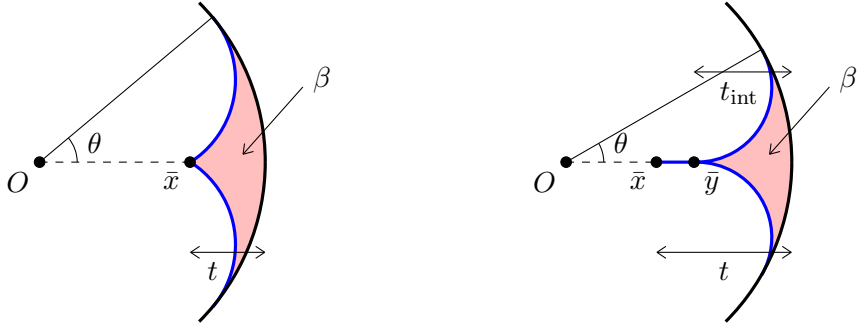

Notice that, in particular, the above characterisation of the minimisers imply that, for every possible $t,\,\beta$, they are unique up to a rotation. Moreover, each minimiser $\bar\gamma^{t,\beta}$ (resp., $\bar\gamma_{t,\beta}$) contains a single point $\bar x$ with $|\bar x|=1+t$ (resp., $|\bar x|=1-t$), and is symmetric with respect to the line $\bar x\R$. In addition, the curve $\bar\gamma^{t,\beta}\setminus \overline B$ (resp., $\bar\gamma_{t,\beta}\setminus \overline B$) meets $\partial B$ in two points, corresponding to an angle $\theta^{t,\beta}$ (resp., $\theta_{t,\beta}$), as in Figure~\ref{FigIsop} and~\ref{Figdentro}. We conclude this step by observing that there exists a purely geometrical constant $C_1>0$ such that, if $t\leq 1$, then
\begin{align}\label{approx1}
\frac 1{C_1}\, t_{\rm ext} \theta^{t,\beta} \leq \beta \leq C_1 t_{\rm ext} \theta^{t,\beta} \,, && \frac 1{C_1} \, t_{\rm int} \theta_{t,\beta}\leq \beta \leq C_1 t_{\rm int} \theta_{t,\beta} \,,
\end{align}
as one can derive by elementary geometrical means recalling that the different parts of $\bar\gamma^{t,\beta}\setminus \partial B$ and $\bar\gamma_{t,\beta}\setminus \partial B$ are arcs of circle.

\step{III}{The inequalities~(\ref{obvapp}).}
In this step, we show that
\begin{align}\label{obvapp}
P(E\cup B) \geq \ell(\bar\gamma^{\delta^+,\nu})\,, && P(E\cap B) \geq \ell(\bar\gamma_{\delta^-,\nu})\,.
\end{align}
First of all, we notice that by~(\ref{defdeltaA}) and~(\ref{precdelta}) the set $E\cup B$ has area $\pi+\nu$ and is contained in the ball $B_{1+\delta^+}$, and similarly $E\cap B$ has area $\pi-\nu$ and contains the ball $B_{1-\delta^-}$. As a consequence, the sets $\Gamma^{\delta^+,\nu}$ and $\Gamma_{\delta^-,\nu}$ are non-empty, so the curves $\bar\gamma^{\delta^+,\nu}$ and $\bar\gamma_{\delta^+,\nu}$ are defined and the inequalities~(\ref{obvapp}) make sense. Moreover, by standard approximation, for every $\sigma>0$ there exists a smooth set $\widetilde E\supseteq B$ with $|P(\widetilde E)-P(E\cup B)|<\sigma$ and such that $|\widetilde E\setminus B|=\nu$ and $\sup\{s>0:\, |\widetilde E\setminus B_s|>0\}=1+\delta^+$. Since $E$ is connected and with no holes in the sense of~(\ref{connectedness}) and~(\ref{imposs}), without loss of generality we can assume that the same is true for $\widetilde E$, hence $\tilde\gamma=\partial\widetilde E$ is a smooth, injective curve. By construction, $\tilde\gamma\in\Gamma^{\delta^+,\nu}$, so
\[
P(E\cup B)+\sigma > P(\widetilde E) = \ell(\tilde\gamma)\geq \ell(\bar\gamma^{\delta^+,\nu})\,,
\]
and since $\sigma$ is arbitrary the left inequality in~(\ref{obvapp}) follows. The proof of the right one is completely similar.

\step{IV}{The curves $\gamma^\pm$.}
In this step, we consider yet another minimisation problem. That is, we minimise $\ell(\gamma)$ among the curves in
\[
\Gamma^+:=\Big\{\gamma\in\Gamma:\, E_\gamma\supseteq B,\, |E_\gamma\setminus B|=\nu,\, \gamma\supseteq \big\{ \eta\in \S^1:\, \theta^{\delta^+,\nu}\leq \eta\leq 2\pi - \theta^{\delta^+,\nu}\big\}\Big\}\,.
\]
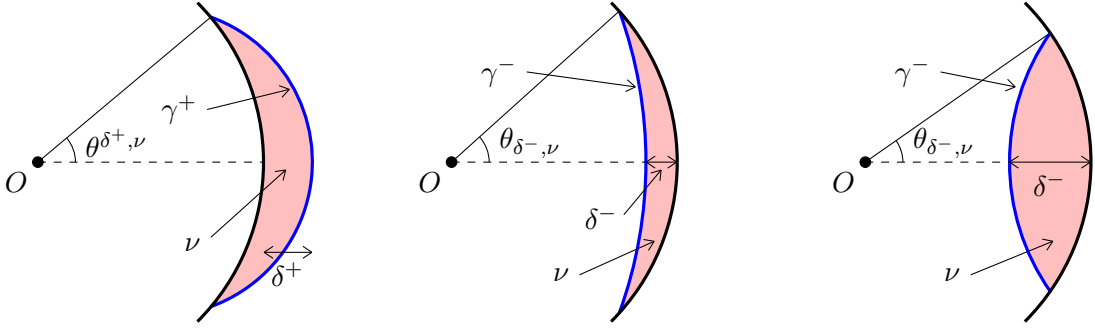
\begin{figure}[htbp]
\begin{tikzpicture}
\fill[fill=red!25!white] (2.3,1.93) arc(70:-70:2.05) arc(-40:40:3);
\filldraw (0,0) circle(2pt);
\draw (0,0) node[anchor=north east] {$O$};
\draw (0,0) -- (2.3,1.93);
\draw (0.5,0) arc(0:40:.5);
\draw (0.5,0.25) node[anchor=west] {$\theta^{\delta^+,\nu}$};
\draw[blue, line width=1.2] (2.3,1.93) arc (70:-70:2.05);
\draw[line width=1.2] (3,0) arc(0:45:3);
\draw[line width=1.2] (3,0) arc(0:-45:3);
\draw[dashed] (0,0) -- (3,0);
\draw[<<<->>>] (3,-1.2) -- (3.65,-1.2);
\draw (3.32,-1.2) node[anchor=north] {$\delta^+$};
\draw (2.3,-1.1) node[anchor=east] {$\nu$};
\draw [->>>] (2.3,-1) -- (3.3,-0.1);
\draw (2.25,.7) node[anchor=east] {$\gamma^+$};
\draw [->>>] (2.25,.7) -- (3.35,.97);
\fill[fill=red!25!white] (7.73,2.01) arc(20:-20:5.87) arc(-42:42:3);
\filldraw (5.5,0) circle(2pt);
\draw (5.5,0) node[anchor=north east] {$O$};
\draw (5.5,0) -- (7.73,2.01);
\draw (6,0) arc(0:42:.5);
\draw (6,0.25) node[anchor=west] {$\theta_{\delta^-,\nu}$};
\draw[blue, line width=1.2] (7.73,2.01) arc (20:-20:5.87);
\draw[line width=1.2] (8.5,0) arc(0:45:3);
\draw[line width=1.2] (8.5,0) arc(0:-45:3);
\draw[dashed] (5.5,0) -- (8.08,0);
\draw[<<<->>>] (8.5,0) -- (8.08,0);
\draw (7.49,-0.45) node[anchor=north] {$\delta^-$};
\draw[->>>] (7.79,-0.7)--(8.29,-.1);
\draw (7.2,-1.5) node[anchor=east] {$\nu$};
\draw [->>>] (7.2,-1.5) -- (8.2,-1.05);
\draw (6.5,1.2) node[anchor=east] {$\gamma^-$};
\draw [->>>] (6.5,1.2) -- (7.95,1);
\fill[fill=red!25!white] (13.46,1.72) arc(145:215:3) arc(-35:35:3);
\filldraw (11,0) circle(2pt);
\draw (11,0) node[anchor=north east] {$O$};
\draw (11,0) -- (13.46,1.72);
\draw (11.5,0) arc(0:35:.5);
\draw (11.5,.25) node[anchor=west] {$\theta_{\delta^-,\nu}$};
\draw[blue, line width=1.2] (13.46,1.72) arc (145:215:3);
\draw[line width=1.2] (14,0) arc(0:45:3);
\draw[line width=1.2] (14,0) arc(0:-45:3);
\draw[dashed] (11,0) -- (12.91,0);
\draw[<<<->>>] (14,0) -- (12.91,0);
\draw (13.45,0) node[anchor=north] {$\delta^-$};
\draw (12.4,-1.5) node[anchor=east] {$\nu$};
\draw [->>>] (12.4,-1.5) -- (13.5,-1.05);
\draw (12,1.2) node[anchor=east] {$\gamma^-$};
\draw [->>>] (12,1.2) -- (13,1);
\end{tikzpicture}
\caption{The minimisers $\gamma^\pm$.}\label{FigYAM}
\end{figure}
By compactness, such a minimiser exists, and we call it $\gamma^+$. Notice that this minimisation problem is trivial, in fact $\gamma^+\setminus \partial B$ is simply the arc of circle which meets $\partial B$ in the two points $(\cos\theta^{\delta^+,\nu},\pm\sin\theta^{\delta^+,\nu})$ and encloses an area $\nu$ outside of $B$, as shown in Figure~\ref{FigYAM}, left. Analogously, we can minimise $\ell(\gamma)$ among curves in
\[
\Gamma^-:=\Big\{\gamma\in\Gamma:\, E_\gamma\subseteq B,\, |B\setminus E_\gamma| = \nu,\, \gamma\supseteq \big\{ \eta\in \S^1:\, \theta_{\delta^-,\nu}\leq \eta\leq 2\pi - \theta_{\delta^-,\nu}\big\}\Big\}\,.
\]
It is again obvious that a minimiser $\gamma^-$ exists, and that $\gamma^-\setminus\partial B$ is the arc of circle which meets $\partial B$ in the two points $(\cos\theta_{\delta^-,\nu},\pm\sin\theta_{\delta^-,\nu})$ in such a way that $B\setminus E_{\gamma^-}$ has area $\nu$. Notice that, depending on $\delta^-$ and $\nu$, the curvature of this arc can be positive (as in Figure~\ref{FigYAM}, centre), null, or negative (as in Figure~\ref{FigYAM}, right). 
We also observe that
\begin{align}\label{wcao}
\ell(\bar\gamma^{\delta^+,\nu}) \geq \ell(\gamma^+)+ 2(\delta^+-\delta^+_{\rm ext})\,, && \ell(\bar\gamma_{\delta^-,\nu}) \geq \ell(\gamma^-)+2(\delta^--\delta^-_{\rm int})\,.
\end{align}
In fact, let us call $\bar\gamma^{\delta^+,\nu}_{\rm short}$ the curve $\bar\gamma^{\delta^+,\nu}\cap \overline{B_{1+t_{\rm ext}}}$. In other words, $\bar\gamma^{\delta^+,\nu}_{\rm short}$ coincides with $\bar\gamma^{\delta^+,\nu}$ in the first three cases of Figure~\ref{FigIsop}, while in the last case $\bar\gamma^{\delta^+,\nu}_{\rm short}$ is obtained by removing the final segment (which is counted twice) from $\bar\gamma^{\delta^+,\nu}$. Similarly, we call $\bar\gamma_{\delta^-,\nu}^{\rm short}$ the curve $\bar\gamma_{\delta^-,\nu}\setminus {B_{1-t_{\rm int}}}$, that is, the whole curve $\bar\gamma_{\delta^-,\nu}$ in the first case of Figure~\ref{Figdentro}, or the curve without the final segment (which is counted twice) otherwise. By construction and by definition of $\theta^{\delta^+,\nu}$ and $\theta_{\delta^-,\nu}$, the curves $\bar\gamma^{\delta^+,\nu}_{\rm short}$ and $\bar\gamma_{\delta^-,\nu}^{\rm short}$ belong to $\Gamma^+$ and $\Gamma^-$ respectively, hence
\[
\ell(\bar\gamma^{\delta^+,\nu})=\ell(\bar\gamma^{\delta^+,\nu}_{\rm short}) + 2(\delta^+-\delta^+_{\rm ext}) \geq \ell(\gamma^+) + 2(\delta^+-\delta^+_{\rm ext})\,,
\]
and similarly
\[
\ell(\bar\gamma_{\delta^-,\nu})=\ell(\bar\gamma_{\delta^-,\nu}^{\rm short}) + 2(\delta^--\delta^-_{\rm int}) \geq \ell(\gamma^-) + 2(\delta^--\delta^-_{\rm int})\,,
\]
so~(\ref{wcao}) is established.\par

Finally, as shown in Figure~\ref{FigYAM}, we will call $d^+$ (resp., $d^-$) the maximal distance between points of $\gamma^+$ (resp., $\gamma^-$) and $\partial B$. As in the previous step, a simple geometric argument ensures the existence of a purely geometric constant $C_2$ such that, if $d^+ \leq 1$, then
\begin{align}\label{approx2}
\frac 1{C_2}\, d^+ \theta^{\delta^+,\nu} \leq \nu \leq C_2 d^+ \theta^{\delta^+,\nu} \,, && \frac 1{C_2} \, d^- \theta_{\delta^-,\nu}\leq \nu \leq C_2 d^- \theta_{\delta^-,\nu} \,.
\end{align}

\step{V}{Conclusion.}
We are now in position to conclude the proof. In fact, let $\bar\delta\ll 1$ be a geometrical constant, to be specified later. Keeping in mind the isoperimetric inequality and~(\ref{centered}), an immediate compactness argument ensures the existence of a constant $\kappa>0$, depending on $\bar\delta$, such that $P(E) \geq 2\pi + \kappa$ if $\delta^++\delta^-\geq \bar\delta$. If $\delta^+\leq 2\pi$ we can then estimate
\[
P(E)-P(B) \geq \kappa \geq \frac \kappa{2\pi+1} \, (\delta^++\delta^-) \geq \frac \kappa{\pi (2\pi+1)} \, \nu (\delta^++\delta^-)\,,
\]
so~(\ref{thesisquad}) is proved if $\delta^++\delta^-\geq \bar\delta$ and $\delta^+\leq 2\pi$. If $\delta^+>2\pi$, then we have $P(E)\geq 2\delta^+$, so that
\[
P(E)-P(B) \geq \delta^+ \geq \frac 1{1+(2\pi)^{-1}}\, (\delta^++\delta^-)\geq \frac 1{\pi+1/2}\, \nu (\delta^++\delta^-)\,,
\]
so once again~(\ref{thesisquad}) is proved.\par

To complete the proof, we have then only to deal with the case when $\delta^++\delta^-\leq \bar\delta$, so from now on we assume this inequality to be true. Notice that we can apply~(\ref{approx1}) with $\beta=\nu$ and with $t=\delta^\pm$, since $\delta^+$ and $\delta^-$ are smaller than $\bar\delta\ll 1$, hence much smaller than $1$. In particular, this implies that $\nu\ll 1$, and it immediately follows that $d^\pm\ll 1$, where $d^\pm$ are given as in Step~IV, so they depend on $\nu,\,\theta^{\delta^+,\nu}$ and $\theta_{\delta^-,\nu}$. As a consequence, we can apply also~(\ref{approx2}), which together with~(\ref{approx1}) gives
\begin{align}\label{dapproxt}
\frac{d^+}{C_1 C_2} \leq \delta^+_{\rm ext} \leq C_1 C_2 d^+ \,, && \frac{d^-}{C_1 C_2} \leq \delta^-_{\rm int} \leq C_1 C_2 d^-\,,
\end{align}
where $\delta^+_{\rm ext}$ and $\delta^-_{\rm int}$ are defined as in Step~II. Let us also call
\begin{align}\label{usefuldeffs}
C_3 = \min\bigg\{ \frac{\sqrt 2}{12C_1C_2},\, \frac 2\pi\bigg\}\,, && \Theta = \min \bigg\{ \frac \pi 4\,,\, \frac{3\pi C_3}{2C_1}\bigg\}\,.
\end{align}
Let us now assume that $\delta^+\geq \delta^-$. If $\theta^{\delta^+,\nu}\geq \Theta$, then we are necessarily in the first case of Figure~\ref{FigIsop}, because otherwise we find a contradiction with the fact that $\nu\ll 1$ (notice that the contradiction holds if the upper bound on $\nu$ is sufficiently small, which in turn is true provided that $\bar\delta$ has been chosen small enough). As a consequence, $\delta^+_{\rm ext}=\delta^+$, and then by the quantitative isoperimetric inequality (see for instance~\cite{FMP,CL}) and~(\ref{approx2}) we get
\[
P(E)-P(B)\geq C_4 \nu^2 \geq \frac{C_4}{C_1}\, \nu \delta^+_{\rm ext} \theta^{\delta^+,\nu}
\geq \frac{C_4\Theta}{C_1}\, \nu \delta^+ 
\geq \frac{C_4\Theta}{2C_1}\, \nu (\delta^+ +\delta^-)
\]
for another geometric constant $C_4$. We have then obtained~(\ref{thesisquad}) and so the proof is concluded if $\delta^+\geq \delta^-$ and $\theta^{\delta^+,\nu}\geq\Theta$. Suppose instead that $\theta^{\delta^+,\nu}<\Theta$, still under the assumption that $\delta^+\geq\delta^-$. In this case, as soon as
\[
\bar\delta < \frac{\sqrt 2}{16 C_1C_2}\,,
\]
by~(\ref{dapproxt}) and~(\ref{usefuldeffs}) we have
\[
d^+ \leq C_1 C_2 \delta^+_{\rm ext} \leq C_1 C_2 \bar\delta <\frac{\sqrt 2}{16}\leq \frac{\cos \Theta}8\leq \frac{\cos(\theta^{\delta^+,\nu})}8\,.
\]
As a consequence, we can apply Lemma~\ref{lempag1} with $\delta=d^+$ and $\bar\theta=\theta^{\delta^+,\nu}$. Observe that with this choice the path $\tau(d^+)$ of Lemma~\ref{lempag1} coincides with $\gamma^+\setminus \partial B$, and $\mu(\delta)=\nu$. Thus, (\ref{obvapp}), (\ref{wcao}), the estimate~(\ref{estpag1}), (\ref{dapproxt}) and~(\ref{usefuldeffs}) give
\[\begin{split}
P(E\cup B) - P(B) &\geq \ell(\bar\gamma^{d^+,\nu}) - 2\pi
\geq \ell(\gamma^+)+2(\delta^+-\delta^+_{\rm ext})-2\pi\\
&=\tau(d^+)- \tau(0) +2(\delta^+-\delta^+_{\rm ext})
\geq \nu + \frac{\cos \bar\theta}6\, d^+\nu +2(\delta^+-\delta^+_{\rm ext})\\
&\geq \nu + \frac{\sqrt 2}{12}\, d^+\nu +\frac{2(\delta^+-\delta^+_{\rm ext})}\pi\,\nu
\geq \nu + \frac{\sqrt 2}{12C_1C_2}\, \delta^+_{\rm ext} \nu +\frac{2(\delta^+-\delta^+_{\rm ext})}\pi\, \nu\\
&\geq \nu + C_3\nu\delta^+\,.
\end{split}\]
On the other hand, by the standard isoperimetric inequality, still minding that $\nu\ll 1$, we have
\[
P(E\cap B) \geq 2\pi \sqrt{\frac{|E\cap B}\pi}
=2\pi \sqrt{\frac{\pi-\nu}\pi}
\geq 2\pi \bigg(1 - \frac \nu{2\pi} - \frac{\nu^2}{6\pi^2}\bigg) = 2\pi - \nu -\frac{\nu^2}{3\pi}\,.
\]
Putting together the last two estimates, also by~(\ref{usefuldeffs}) we have then
\[
P(E)-P(B) \geq P(E\cup B) + P(E\cap B) -2P(B)
\geq C_3\nu\delta^+ -\frac{\nu^2}{3\pi}
\geq \frac{C_3}2\, \nu\delta^+ \geq \frac{C_3}4\, \nu(\delta^++\delta^-)\,.
\]
We have then obtained~(\ref{thesisquad}) under the assumption that $\theta^{\delta^+,\nu}>\Theta$ and $\delta^+\geq\delta^-$, so the proof is completed for the case $\delta^+\geq\delta^-$.\par

The proof for the case $\delta^-\geq\delta^+$ is exactly the same, just replacing in the obvious way $\delta^+,\,\delta^+_{\rm ext}$ and $\theta^{\delta^+,\nu}$ with $\delta^-,\,\delta^-_{\rm int}$ and $\theta_{\delta^-,\nu}$.
\end{proof}

We are finally in position to prove Theorem~\ref{2dcase}.

\proofof{Theorem~\ref{2dcase}}
Let us fix a small, positive constant $\bar\eps$, and let $0<\eps<\bar\eps$. Let also $E$ be a minimiser of the functional $\widetilde\F_\eps$, which exists by Proposition~\ref{genex}. Provided that $\bar\eps<\eps_2$, we can apply Lemma~\ref{solconn} and Lemma~\ref{noholes}, to find that $E$ is also a minimiser of $\F_\eps$, and that it is connected and with no holes. Up to a translation, we can assume that~(\ref{centered}) holds. It is then possible to apply Lemma~\ref{AdelC}, hence~(\ref{thesisquad}) gives
\begin{equation}\label{2ndthesisquad}
P(E)-P(B) \geq \frac 1C\, |E\Delta B|(\delta^++\delta^-)\,,
\end{equation}
where $C$ is a purely geometric constant and $\delta^\pm$ are defined as in~(\ref{precdelta}). We call
\begin{align*}
E^+=E\setminus B\,, && E^-=B\setminus E\,,
\end{align*}
so that $E=B\cup E^+\setminus E^-$, and since $g$ is positive definite in the sense of Definition~\ref{defpos} we can evaluate
\begin{equation}\label{firstpart}\begin{split}
\En(E)&=\En(E,B)+\En(E,E^+)-\En(E,E^-)\\
&=\En(B)+\En(B,E^+)-\En(B,E^-)+\En(E,E^+)-\En(E,E^-)\\
&=\En(B)+2\big(\En(B,E^+)-\En(B,E^-)\big)+\En(E^+,E^+)-2\En(E^+,E^-)+\En(E^-,E^-)\\
&\geq \En(B)+2\big(\En(B,E^+)-\En(B,E^-)\big)\,.
\end{split}\end{equation}
Notice now that $\delta^+$ and $\delta^-$ are bounded. Indeed, $\delta^-\leq 1$, and clearly $P(E)\geq 2(\delta^+-1)$, so that $\delta^+$ is bounded because Lemma~\ref{lemmafraenkel} implies that $P(E)$ is close to $2\pi$ if $\bar\eps$ is small enough. Moreover, since the assumption~(\ref{ipLip}) coincides with~(\ref{N-2}) because $N=2$, Lemma~\ref{lemadm} implies that the function $\Phi:\R^+\to\R^+$ defined as in~(\ref{defPhi}) is locally Lipschitz continuous. Therefore, there exists a constant $L\in\R^+$ such that
\begin{align*}
\Phi(|x|)\geq \Phi(1)-L(|x|-1) \quad \forall\, x\in E^+\,, &&
\Phi(|x|)\leq \Phi(1)+L(1-|x|) \quad \forall\, x\in E^-\,.
\end{align*}
Keeping in mind that $|E^+|=|E^-|=|E\Delta B|/2$, we can then evaluate
\[
\En(B,E^+) = \int_{E^+} \Phi(|x|)\,dx \geq \Phi(1)|E^+| - L\int_{E^+} |x|-1\,dx
\geq \Phi(1)|E^+| - L\delta^+|E^+|\,,
\]
and similarly
\[
\En(B,E^-) = \int_{E^-} \Phi(|x|)\,dx \leq \Phi(1)|E^-| + L\int_{E^-} 1-|x|\,dx
\leq \Phi(1)|E^+| + L\delta^-|E^-|\,,
\]
so that~(\ref{firstpart}) gives
\[
\En(E)\geq \En(B)+2\big(\En(B,E^+)-\En(B,E^-)\big)
\geq \En(B)-L|E\Delta B|(\delta^++\delta^-)\,.
\]
Putting this estimate together with~(\ref{2ndthesisquad}), and keeping in mind that $E$ is a minimiser of $\F_\eps$, we find then
\[
\F_\eps(B)\geq \F_\eps(E)=P(E) + \eps\En(E)\geq \F_\eps(B) +\bigg(\frac 1C-\eps L\bigg)\, |E\Delta B|(\delta^++\delta^-)\,.
\]
As a consequence, provided that $\bar\eps<(CL)^{-1}$, we obtain that $E$ coincides with $B$. The proof is then concluded.
\end{proof}


\end{document}